\documentclass[reqno]{amsart}
\usepackage{amsmath}
\usepackage{amssymb}
\usepackage{stmaryrd}
\title{
Removal Lemma for Infinitely-Many Forbidden Hypergraphs
 and Property Testing
}
\author{Yoshiyasu Ishigami}
\address{Department of Information and Communication Engineering, 
The University of Electro-Communications, Chofu, Tokyo 182-8585, Japan.
}
\email{yoshiyas@ice.uec.ac.jp}
\date{\today}
\def\picture #1 by #2 (#3){
                \vbox to #2{
                        \hrule width #1  height 0pt depth
0pt
                        \vfill
                        \special{picture #3}}}
\def\scaledpicture #1 by #2(#3 scaled #4){{
                \dimen0=#1 \dimen1=#2
                \divide\dimen0 by 1000 \multiply\dimen0 by
#4
                \divide\dimen1 by 1000 \multiply\dimen1 by
#4
                \picture\dimen0 by \dimen1 (#3 scaled #4)}}
\newtheorem{adf}{Definition}[section]
\newtheorem{tha}{Theorem}[section]

\newtheorem{thm}{Theorem}[section]
\newtheorem{cor}[thm]{Corollary}
\newtheorem{lemma}[thm]{Lemma}
\newtheorem{remark0}[thm]{Remark}

\newtheorem{algorithm0}[thm]{Algorithm}
\newtheorem{asett}[thm]{Setup}

\newenvironment{df}{
\begin{adf}\begin{sl}}
{\end{sl}\thqed\end{adf}}
\newenvironment{sett}{
\begin{asett}\begin{sl}}{\end{sl}\thqed\end{asett}}

\newcommand{\thqed}{\hfill\fbox{}\\ }
\newcommand{\Proof }{{\bf Proof} : }

\newcommand{\lemmaqed}{\hfill\bigskip\fbox{}\\}

\def\cil#1{\left\lceil{#1}\right\rceil}
\def\brkt#1{\left({#1}\right)}

\def\ang#1{\langle{#1}\rangle}

\newcommand{\naturalset}{{\mathbb N}}

\newcommand{\Prob}{{\mathbb P}}
\newcommand{\Ex}{{\mathbb E}}

\newtheorem{procedure0}{Procedure}

\newcount\timehh\newcount\timemm\timehh=\time\divide\timehh by 60
\timemm=\time\count255=\timehh\multiply\count255 by-60
\advance\timemm by \count255
\ifnum\timehh=12\def\apm{pm}\else
\ifnum\timehh>12\def\apm{pm}\advance\timehh by-12\else
\def\apm{am}\fi\fi
\def\timestamp{\number\timehh\,:\,\ifnum\timemm<10 0\fi\number\timemm\,\apm}
\begin{document}
\maketitle
\begin{abstract}
We prove a removal lemma for infinitely-many forbidden hypergraphs.
It affirmatively settles 
a question on property testing raised by Alon and Shapira (2005) \cite{AS,AS05}.
All monotone hypergraph properties and all hereditary partite hypergraph properties 
are testable.
Our proof constructs a constant-time probabilistic algorithm to edit a small number of edges.
It also gives a quantitative 
bound in terms of a coloring number of the property.
It is based on a new hypergraph regularity lemma \cite{I06}.
\end{abstract}
\section{Introduction}
The research of removal lemmas has started by Rusza and Szemer\'edi \cite{RuSz}, who 
considered an ordinary graph case.
Frankl and R\"{o}dl \cite{FR02} suggested that if a hypergraph version of 
the removal lemma can be proven, it 
yields Szemer\'edi's famous theorem on arithmetic progressions 
\cite{Sz75}. They actually showed an alternative
 proof of the theorem for length four \cite{Sz69} by 
showing a $3$-uniform hypergraph regularity lemma. 
Later Solymosi \cite{So03,So04} showed that the $k$-uniform hypergraph removal lemma 
(conjecture)  implies 
not only Szemer\'edi's theorem but also its multidimensional extension by 
Furstenberg and Katznelson \cite{FK78}, which had been proved only by ergodic theory until 
recently. 
Finally Gowers \cite{G} and R\"{o}dl and his collaborators \cite{RSk04,NRS} 
 obtained  the hypergraph removal lemma as a corollary of their regularity lemmas. 
Slightly later, Tao \cite{Tao06} gave another regularity lemma,
 which yields the hypergraph removal lemma.
Very recently \cite{I06} gave a new regularity lemma with a clear 
construction and a simple proof, which we will use in this paper.
\par
Tao \cite{Tao06b} gave another proof of hypergraph removal lemma by using 
ergodic theoretic ideas. It is 
nonconstructive but is independent from any regularity lemma.
\par
These hypergraph removal lemmas deal with one forbidden hypergraph.
It is straightforward to rewrite them for a finite family of forbidden hypergraphs.
For details of hypergraph removal lemmas, see \cite[\S11.6.pp.454-463]{TV}.
\par
In the below, a partite hypergraph is {\bf $h$-vertex} if and only if 
each partite set contains exactly $h$ vertices.
The main puropose of this paper is to show the following.
\begin{thm}[Removal lemma for infinitely-many forbidden $r$-partite hypergraphs]
\label{6S15g}
Let $r\ge k$ be positive integers and $\varepsilon>0$.
Let ${\mathcal F}=\bigcup_{h=1}^\infty {\mathcal F}_h$ where 
${\mathcal F}_h$ is a 
family of $h$-vertex $k$-uniform (${r}$-partite hyper)graphs.
Then there exist 
constants $c>0$ and $h_0$ 
such that for any integer $N$ if 
 ${\bf G}$ is an $N$-vertex $k$-uniform (${r}$-partite hyper)graph then 
 at least one of the following two holds.
\\
{\rm (i)} One can modify at most $\varepsilon {r\choose k}N^k$ edges 
so that the hypergraph obtained from ${\bf G}$ does not have 
a copy of any member of ${\mathcal F}$ as an induced-subhypergraph.
\\
{\rm (ii) } There exist $h\le h_0, F\in {\mathcal F}_h$ such that 
${\bf G}$ contains at least $cN^{rh}$ copies of $F$.
\end{thm}
Our main theorem, Theorem \ref{6S05}, will be presented in a more general frame.
For example, (1) each forbidden hypergraph $F$ may contain not only black and white 
edges but also \lq invisible\rq\ edges, (so it contains the both cases of 
indueced and non-induced subgraphs)
(2)  $F$ and ${\bf G}$ may have a constant 
number of colors other than black and white, and (3) partite vertex sets of ${\bf G}$ 
will be discrete probability spaces with finite vertices, where the sizes of 
two partite vertex sets may not be equal. 
Those are not insiginificant extensions. We employ them to make our 
argument natural.
\section{Testability}
Property testing was firstly considered by Blum et al.\cite{BLR}, and 
the general notion of property testing was first given 
by Rubinfeld and Sudan in \cite{RuSu}.
Goldreich et al. \cite{GGR} firstly investigated it in combinatorial context, in which 
they considered ordinary graphs.
See surveys \cite{Ron,Fis,AS}.
\begin{df}[Hypergraph property]
Two hypergraphs are {\bf isomorphic} if and only if 
one can be equivalent to the other by some bijection (permutation) 
between the two vertex sets.
(For the case of $r$-partite hypergraphs, the bijection should be 
\lq partitionwise\rq, i.e. any vertex in any partite set cannot be 
replaced to a different partite set, and furthermore 
partite sets have their own labels from $1,\cdots,r$,  any of which we cannot change.
)
A {\bf hypergraph property} (or property, simply) is a 
class of hypergraphs such that if a hypergraph belongs to the class 
(satisfies the property) then 
any other hyprgraphs isomorphic to the hypergraph belong to it.
\par
A hypergraph property is {\bf monotone} if and only if 
when a hypergraph satisfies the property, any (induced/non-induced) subgraph of it 
satisfies the property.
A hypergraph property is {\bf hereditary} if and only if 
when a hypergraph satisfies the property, any induced subgraph of it 
satisfies the property. 
In other words, a monotone (or hereditary) property is colosed under 
any deletion of vertices and edges (or vertices, respectively).
Clearly any monotone property is hereditary.
\end{df}
\begin{df}[$\epsilon$-far]
A $k$-uniform hypergraph is 
{\bf $\epsilon$-far} from a property ${\mathcal P}$ if and only if 
the hypergraph cannot satisfy ${\mathcal P}$ even after modifying at most 
$\epsilon$ portion of edges of the underlying complete hypergraph.
(i.e. it is $\epsilon {N\choose k}$ edges for $N$-vertex hypergraphs 
(with no vertex partitions) and 
is $\epsilon {r\choose k}N^k$ edges for $r$-partite hypergraphs 
with $N$ vertices in each partite set.)
\end{df}
\begin{df}[Property test]
A property is {\bf testable} if and only if 
there exists a randomized algorithm such that, for any $\epsilon>0$ and 
any object(a hypergraph) given as inputs,
if \\
(1) the input object satisfies the property
 or \\
(2) it is $\epsilon$-far from the property\\
then 
with probability at least 0.9 the algorithm correctly answers 
which case of the two it is, in a constant time independent from the size 
of the object (the number of vertices in the hypergraph).
(The time can depend on $\epsilon.$)
A testable property is {\bf testable with one-sided error} if and only if 
its answer is correct always(with probability 1)
 whenever the input satisfies the property (i.e. the case (1)).
\end{df}
\begin{thm}\label{6S15h}
Every hereditary property of constant-partite hypergraphs is testable with 
one-sided error.
\end{thm}
\begin{proof}
{\bf [Design of the algorithm]} 
Firstly we will present a random algorithm for hereditary property ${\mathcal P}$ 
with one-sided error.
Fix $r\ge k$ and $\epsilon>0$.
Let ${\mathcal F}_h$ be the set of all $h$-vertex 
$r$-partite $k$-uniform hypergraphs 
which do not satisfy ${\mathcal P}$. Let ${\mathcal F}=\bigcup_{h\ge 1}{\mathcal F}_h$.
With these parameters, Theorem \ref{6S15g} gives us constants $c>0$ and $h_0$.
Our algorithm goes as follows.
Given the input hypergraph ${\bf G}$, the algorithm randomly chooses 
vertices $W^{(i)}=\dot{\bigcup}_{j\in [r]}\{v^{(i)}_{j,j'}\in 
{\bf \Omega}_j| j'\in [h_0]\}$ 
for $\cil{3/c}$ times $(i\in [\cil{3/c}])$,
where ${\bf \Omega}_j$ denotes the $j$-th partite vertex-set of ${\bf G}$.
Then it declares ${\bf G}$ to satisfy ${\mathcal P}$ if-and-only-if, for all $i,$
 the $h_0$-vertex hypergraph $H^{(i)}$induced by $W^{(i)}$ satisfies 
${\mathcal P}$(i.e. it is not isomorphic to any member of ${\mathcal F}_{h_0}$).
\par
{\bf [Verification of the algorithm]}
Suppose that the input ${\bf G}$ satisfies ${\mathcal P}$.
Since ${\mathcal P}$ is hereditary, 
all induced sub(hyper)graphs of ${\bf G}$ (thus also all $H^{(i)}$) satisfy ${\mathcal P}.$
Thus the algorithm declares correctly with probability one.
\par
Assume that ${\bf G}$ is $\epsilon$-far from ${\mathcal P}$.
Theorem \ref{6S15g} says that there exists an $F\in {\mathcal F}_{h}$ 
with an $h\le h_0$ such that ${\bf G}$ contains at least $cN^{rh}$ copies of $F.$
Let $F^+$ be the $h_0$-vertex hypergraph obtained from $F$ by adding some 
isolated vertices. Since all added vertices are isolated, 
${\bf G}$ contains at least $cN^{rh_0}$ copies of $F^+.$
Consequently for any fixed $i,$ the probability that 
$H^{(i)}$ is isomorphic to $F^+$ is 
at least $c.$
\par
 If $F^+\not\in {\mathcal F}_{h_0}$ then $F^+$ satisfies 
${\mathcal P}$ and then its induced-subgraph 
$F$ also satisfies ${\mathcal P}$, contradicting $F\in {\mathcal F}_h.$
Thus $F^+\in {\mathcal F}_{h_0}$ and $F^+$ does not satisfy ${\mathcal P}$.
If some $H^{(i)}$ is isomorphic to $F^+$ 
then $H^{(i)}$ does not satisfy ${\mathcal P}$ and 
the algorithm must say that ${\bf G}$ does not satisfy ${\mathcal P}$.
Thus the probability that the algorithm outputs the wrong answer is 
\begin{eqnarray*}
\Prob[\mbox{The algorithm incorrectly says that } {\bf G} \mbox{ is in } {\mathcal P}]
&\le & \Prob[H^{(i)}\mbox{ is not isomorphic to }F^+ \mbox{ for all } i]\\
&\le & (1-c)^{3/c}\\
&\le & e^{-3}\\
&<& 0.1.
\end{eqnarray*}
\end{proof}
Alon and Shapira \cite{AS05} showed that every monotone graph property is 
testable with one-sided error, where their proof using graph regularity lemma 
\cite{Sz} was (probabilistically-)constructive and gave a quantitative bound. 
Lov\'asz and Szegedy \cite{LS05} gave an alternative proof by using 
graph sequences \cite{LS04}, which is short but not constructive.
(Thus the input graph is not $\epsilon$-far but does not satisfy the target property, 
they do not give us any procedure about the way of modifying the input graph so that 
the resulting graph satisfies the target property.
)
Alon and Shapira \cite{AS,AS05}
asked whether it can be extended to 
uniform hypergraphs \cite{AS,AS05}.
Their main interest seems to be 
whether recently discovered hypergraph regularity lemmas 
(\cite{RSk04,NRS,G,Tao06,RS}) are strong enough for applications to property testing.
It had been known that they are strong enough for Szemer\'edi theorem on 
arithmetic progressions (\cite{Sz75}) and its variants.
Avart et al. \cite{ARS} showed it for 3-uniform hypergraphs, by developing 
their argument with the 3-uniform hypergraph regularity lemma of \cite{FR02}.
We will answer their question as follows, by using a new hypergraph regularity platform 
\cite{I06}.
\begin{cor}
Every monotone property of hypergraphs with no vertex partitions 
is testable with one-sided error.
\end{cor}
\begin{proof}
 It easily follows from Theorem \ref{6S15h}. 
Choose an $r$ to be a constant large enough with respect to 
$1/\epsilon$ and $k$. We modify the input non-partite hypergraph to be an $r$-partite 
hypergraph by decomposing the vertex set to $r$ disjoint vertex partite-sets and 
by deleting(invisualizing) \lq non-partitionwise\rq\ edges 
(i.e. deleting any edge with at least two vertices being in a common vertex partite-set).
Note that there are at most $r\cdot r^{k-2}N^k<0.1 \epsilon {r\choose k}N^k$ such deleted 
edges.
It is reduced to Theorem \ref{6S15h}.
\end{proof}
In the course of writing the first draft of this paper \cite{I06m}, I learned that 
R\"{o}dl and Schacht \cite{RSgene,RSgene2} proved the above independently from me.
Their method even yields that every hereditary non-partite hypergraph properties 
are testable 
with one-sided error. In this sense, their result is stronger.
However the approaches are siginificantly different. 
They combined their regularity lemma with 
 the (non-constructive) idea of graph limits from \cite{LS04,LS05}, 
without extending the approaches of \cite{AS05,ARS}.
It may be practically impossible or hard to show the testability even 
for monotone hypergraph properties 
by extending the proof of \cite{AS05, ARS} naturally under their regularity lemma.
Their proof is not constructive and 
does not give any quantitative bound.
On the other hand, our proof gives a procedure about which edges should be modified 
in the given hypergraph. A quantitative bound for the number of edges 
to modify  can be calculated in terms of a {\em coloring number} of the property, 
though the bound seems to be weak. 
Improving the bound would be an interesting research theme.
Their proof is based on their heavy regularity lemma, while 
ours is based on a new regularity lemma \cite{I06}, which 
has a shorter proof.
\section{Statement of the Main Theorem}
In this paper, we denote by $\Prob$ and $\Ex$ the probability and expectation, 
respectively. We denote the conditional probability and exepctation by
$\Prob[\cdots|\cdots]$ and $\Ex[\cdots|\cdots].$ 
\begin{sett}\label{mt051231}
Throughout this paper, we fix a positive integer $r$ and 
an \lq index\rq\ set $\mathfrak{r}$ with $|\mathfrak{r}|=r.$ 
Also we fix a probability space 
$({\bf \Omega}_i,{\mathcal B}_i,\Prob)$ 
for each $i\in \mathfrak{r}$.
Assume that ${\bf \Omega}_i$ is finite (but its cardinality may not be 
constant) 
and ${\mathcal B}_i=2^{{\bf \Omega}_i}$ 
for the sake of simplicity.
Write ${\bf \Omega}:=({\bf \Omega}_i)_{i\in \mathfrak{r}}$.
\end{sett}
In order to avoid using measure-theoretic jargons like measurability or 
Fubini's theorem frequently to readers who are interested only in applications to 
discrete mathematics, 
we assume ${\bf \Omega}_i$ as a (non-empty) finite set.
However our argument should be extendable to a more general probability space.
For applications, ${\bf \Omega}_i$ would contain a huge number of vertices, though 
we do not use the assumption in our proof.
\par
For an integer $a$, we write $[a]:=\{1,2,\cdots,a\},$ and 
${\mathfrak{r}\choose [a]}:=\dot{\bigcup}_{i\in [a]}{\mathfrak{r}\choose i}
=\dot{\bigcup}_{i\in [a]}\{I\subset \mathfrak{r}| |I|=i\}.$
When $r$ sets $X_i, i\in {\mathfrak r},$ with indices from ${\mathfrak r}$ are 
called {\bf vertex sets}, 
we write $X_J:=\{e\subset \dot{\bigcup}_{i\in J}X_i| |e\cap X_j|=1 \forall j\in J\}$ 
whenever $J\subset {\mathfrak r}$.
\begin{df}[(Colored hyper)graphs]
Suppose Setting \ref{mt051231}. 
A {\bf $k$-bound $(b_i)_{i\in [k]}$-colored ($\mathfrak{r}$-partite hyper)graph} $H$ 
is a triple $((X_i)_{i\in\mathfrak{r}},({C}_I)_{I\in {\mathfrak{r}\choose [k]}},
(\gamma_I)_{I\in {\mathfrak{r}\choose [k]}}
)$ where (1) each $X_i$ is a set called a \lq vertex set,\rq\ (2)
${C}_I$ is a set with at most $b_{|I|}$ elements, and 
(3) $\gamma_I$ is a map from $X_I$ to ${C}_I.$
We write $V(H)=\dot{\bigcup}_{i\in \mathfrak{r}}X_i
$ and 
${\rm C}_I(H)={C}_I.$ 
Each element of $V(H)$ is called a {\bf vertex}.
Each element $e\in V_I(H)=X_I, I\in {\mathfrak{r}\choose [k]},$ is called 
 an {\bf (index-$I$) edge}.
Each member in ${\rm C}_I(H)$ is a {\bf (face-)color (of index $I$)}.
Write $H(e)=\gamma_I(e)$ for each $I.$ 
\par
Let $I\in {\mathfrak{r}\choose [k]}$ 
and $e\in V_I(H).$ 
For another index $\emptyset\not=
J\subset I$, we denote by $e|_J$ the index-$J$ edge 
$e\setminus \brkt{\bigcup_{j\in I\setminus J}X_j}\in V_J(H)$.
We define the {\bf frame-color}
and {\bf total-color} of $e$ by ${H}(\partial
e):=({H}(e|_J)|\,\emptyset\not=J\subsetneq I)$ and by
${H}(\ang{e})=H\ang{e}:=({H}(e|_J)|\,\emptyset\not=J\subset I).$
Write ${\rm TC}_I(H):=\{H\ang{e}|{e}\in X_I\}$ and 
$\partial {\rm C}_I(H):=\{H(\partial e)|{e}\in X_I\}.$ 
\par
A {\bf ($k$-bound) simplicial-complex} is a $k$-bound
(colored ${\mathfrak r}$-partite hyper)graph 
such that for each $I\in {\mathfrak{r}\choose [k]}$ 
there exists at most one index-$I$ color called \lq invisible\rq\ 
and that if (the color of) an edge $e$ is invisible then 
any edge $e^*\supset e$ is invisible. An edge or its color 
is {\bf visible} if it is not invisible.
\par
For a $k$-bound graph ${\bf G}$ on ${\bf \Omega}$ and $s\le k$, 
let ${\mathcal S}_{s,h,{\bf G}}$ be the set of $s$-bound 
simplicial-complexes $S$ such that 
(1) each of $r$ vertex sets contains exactly $h$ vertices and 
that (2)
for any $I\in {\mathfrak{r}\choose [s]}$ 
 there is an injection from the index-$I$ visible colors of $S$ to the 
index-$I$ colors of ${\bf G}$.
(When a visible color $\mathfrak{c}$ of $S$ corresponds to another color $\mathfrak{c}'$ of 
${\bf G}$, we simply write $\mathfrak{c}=\mathfrak{c}'$ 
without presenting the injection explicitly.)
For $S\in {\mathcal S}_{s,h,{\bf G}}$, we denote by ${\mathbb V}_I(S)$ the set of 
index-$I$ visible edges. Write ${\mathbb V}_i(S):=\bigcup_{I\in {\mathfrak{r}\choose 
i}}{\mathbb V}_I(S)$ and ${\mathbb V}(S):=\bigcup_i {\mathbb V}_i(S).$
\end{df}
Informally speaking, our aim will be to embed a \lq child graph\rq\ $S$ to 
a \lq mother graph\rq\ ${\bf G}$ on vertex set ${\bf \Omega}$. 
We will use bold fonts for vertices and edges of the mother graph.
\begin{df}[Partitionwise maps]
A {\bf partitionwise map} $\varphi$ is 
a map 
 from $r$ vertex sets $W_i,i\in \mathfrak{r},$ with $|W_i|<\infty$ to {\em the} 
$r$ vertex sets (probability spaces) $
{\bf \Omega}_i,i\in\mathfrak{r}
$,
such that 
each $w\in W_i$ is mappped into ${\bf \Omega}_i.$
We denote by $\Phi((W_i)_{i\in\mathfrak{r}})$ or $\Phi(\bigcup_{i\in\mathfrak{r}}W_i)$ 
the set of partitionwise maps from $(W_i)_i.$
When $W_i=[h]$ or $W_i$ are obvious and $|W_i|=h$, we denote it by 
$\Phi(h)$.
A partitionwise map is {\bf random} if and only if 
each vertex $w\in W_i$ is mutually-independently mapped at random according to 
 the probability space 
${\bf \Omega}_i$.
\par
Define $\Phi(m_1,\cdots,m_{k-1}):=\Phi(m_1)\times\cdots\times\Phi(m_{k-1}).$ 
\end{df}
\begin{df}[$k$-uniform graphs]
A {\bf $k$-uniform $b_k$-colored ($\mathfrak{r}$-partite hyper)graph} is 
a $k$-bound $(b_i)_{i\in [k]}$-colored graph such that \\
(1) 
 if $i<k$ then $b_i=1$ and the unique color is called invisible 
and \\
(2) for each $I$ with $|I|=k$, there is at most one
 index-$I$ color which is called invisible. 
(Note that this word \lq invisible\rq\ 
is slightly different from the same word used in 
the definition of simplicial-complexes.)
Denote by $\mathbb{V}(F)$ the set of 
visible edges of a $k$-uniform graph $F$, 
where a visible edge means an edge whose color is not invisible.
It is called {\bf $h$-vertex} if each partite set contains exactly $h$ vertices.
\end{df}
\begin{thm}[Main Theorem]
\label{6S05}
Let $r\ge k$ and $\vec{b}=(b_i)_{i\in [k]}$ be positive integers and $\varepsilon>0$.
Let ${\mathcal F}=\bigcup_{h=1}^\infty {\mathcal F}_h$ 
 where ${\mathcal F}_h$ is a  
family of $h$-vertex $k$-uniform $b_k$-colored ($\mathfrak{r}$-partite hyper)graphs.
Then there exist 
constants $c>0$ and $h_0$ 
with the following.\par
Let ${\bf G}$ be a $k$-bound $\vec{b}$-colored ($\mathfrak{r}$-partite hyper)graph on 
${\bf \Omega}=({\bf \Omega}_i)_{i\in \mathfrak{r}}.$
Then at least one of the following two holds.
\\
{\rm (i)}
 There exists a $k$-bound $\vec{b}$-colored 
($\mathfrak{r}$-partite hyper)graph ${\bf G'}$ on 
${\bf \Omega}$ such that 
\begin{eqnarray*}
\Prob_{{\bf e}\in {\bf \Omega}_I
}[{\bf G'}({\bf e})\not={\bf G}({\bf e})]\le \varepsilon \quad 
\forall I\in {\mathfrak{r}\choose k}
\end{eqnarray*}
and that for all $h, F\in {\mathcal F}_h,$
\begin{eqnarray*}
\Prob_{\phi\in\Phi(h)}[{\bf G'}(\phi(e))=F(e)\forall e\in {\mathbb V}(F)]=0.
\end{eqnarray*}
{\rm (ii) } There exist $h\le h_0, F\in {\mathcal F}_h$ such that 
\begin{eqnarray*}
 \Prob_{\phi\in\Phi(h)}[{\bf G}(\phi(e))=F(e)
\forall e\in {\mathbb V}(F)]\ge c.
\end{eqnarray*}
\end{thm}
Our proof is constructive. In hypergraph regularity setup of \cite{I06}, 
we will develop the argument which Alon and Shapira \cite{AS05} used 
for graphs. 
\section{Definitions of Regularities and Statement of Regularity Lemma}
\begin{df}[Regularization]
Let $m\ge 0$ and $\varphi
\in \Phi(m).$ 
Let ${\bf G}$ be a $k$-bound graph on ${\bf \Omega}.$
For an integer $1\le s< k$,
the {\bf $s$-regularization} ${\bf G}/^s\varphi$ is
the $k$-bound graph on ${\bf \Omega}
$ obtained from ${\bf G}$ by redefining the color of 
each edge ${\bf e}\in {\bf \Omega}_I$ with $I\in {\mathfrak{r}\choose [s]}$
by the vector
\begin{eqnarray*}
\brkt{{\bf G}/^s\varphi}({\bf e}):=({\bf G}({\bf e}\dot{\cup}{\bf f})
|
J\in {\mathfrak{r}\setminus I\choose [0,s+1-|I|]},
{\bf f}\in {\bf \Omega}_J \mbox{\rm with } 
{\bf f}\subset \varphi(\mathbb{D})
).
\end{eqnarray*}
In the above, when $J=\emptyset,$ we assume ${\bf f}=\emptyset.$
(The sets of colors are naturally extended
while any edge containing at least $s+1$ vertices does not change its (face-)color.)
\par
When $s=k-1$, we simply write ${\bf G}/\varphi={\bf G}/^{k-1}\varphi.$
\par
For $\vec{\varphi}=(\varphi_i)_{i\in [k-1]}\in\Phi(m_1,\cdots,m_{k-1}),$
we define the {\bf regularization of ${\bf G}$ by $\vec{\varphi}$} by
\begin{eqnarray*}
{\bf G}/\vec{\varphi}:=
(({\bf G}/^{k-1}\varphi_{k-1})/^{k-2}\varphi_{k-2})\cdots/^1\varphi_1.
\end{eqnarray*}
\end{df}
\begin{df}[Regularity]
Let ${\bf G}$ be a $k$-bound graph on $
{\bf \Omega}$.
For 
$\vec{\mathfrak{c}}=(\mathfrak{c}_J)_{J\subset I}\in {\rm TC}_I({\bf G}), 
I\in {\mathfrak{r}\choose [k]}$, we define {\bf relative density}
\begin{eqnarray*}
{\bf d}_{\bf G}(\vec{\mathfrak{c}})=
{\bf d}_{\bf G}(\mathfrak{c}_I|(\mathfrak{c}_J)_{J\subsetneq I}):=
\Prob_{{\bf e}\in {\bf \Omega}_I
}[
{\bf G}({\bf e})=\mathfrak{c}_I
|
{\bf G}(\partial{\bf e})=
(\mathfrak{c}_J)_{J\subsetneq I}
].
\end{eqnarray*}
\par
For a nonnegative integer $h$ and ${\epsilon}\ge 0$, 
we say that ${\bf G}$ is 
 {\bf $({\epsilon},k,h)$-regular} 
(or {\bf $({\epsilon},h)$-regular}) if and only if 
there exists a function
${ \delta}: 
{\rm TC}({\bf G})\to [0,\infty)$ 
such that 
\begin{eqnarray*}
\hspace{-5mm}
{\rm (i)}&
 \Prob_{\phi\in\Phi(h)}
[{\bf G}(\phi(e))=S(e)
\forall e\in {\mathbb V}(S)]=
\displaystyle\prod_{e\in {\mathbb V}(S)}
\brkt{
{\bf d}_{\bf G}(S\ang{e})
\dot{\pm}
\delta(S\ang{e})
}&
\forall S\in {\mathcal S}_{k,h,{\bf G}},
\label{mt060719}
\\
\hspace{-5mm}
{\rm (ii)}& \Ex_{{\bf e}\in {\bf \Omega}_I}[\delta({\bf G}\ang{\bf e})]\le 
\epsilon /|{\rm C}_I({\bf G})|
& \forall I\in {\mathfrak{r}\choose [k]},
\label{mt060726}
\end{eqnarray*}
where $a\dot{\pm}b$ means (the interval of) numbers $c$ with 
$\max\{0,a-b\}\le c\le \min\{1,a+b\}$.
Denote by ${\bf reg}_{k,h}({\bf G})$ the minimum value 
of $\epsilon$ such that ${\bf G}$ is $(\epsilon,k,h)$-regular.
\par
For nonnegative integers $h,L$ and $\epsilon\ge 0$, we say that 
 ${\bf G}$ is 
 {\bf $({\epsilon},k,h,L)$-regular} (or {\bf $({\epsilon},h,L)$-regular})
if and only if  ${\bf G}$ is 
 $(\epsilon,k,h)$-regular and the following holds for all $I\in 
{\mathfrak{r}\choose [k]}$ :
\begin{eqnarray*}
\Ex_{\varphi\in\Phi(L)}
\Ex_{{\bf e^*}\in {\bf \Omega}_I}
[
\sum_{\mathfrak{c}_I\in {\rm C}_I({\bf G})
}
\brkt{
{\bf d}_{{\bf G}/\varphi}(\mathfrak{c}_I|{\bf G}/\varphi(\partial {\bf e^*}))
-
{\bf d}_{\bf G}(\mathfrak{c}_I|{\bf G}(\partial {\bf e^*}))
}^2
]\le \brkt{
{\epsilon\over |{\rm C}_I({\bf G})|}
}^2
\end{eqnarray*}
where we naturally write 
\begin{eqnarray}
{\bf d}_{{\bf G}/\varphi}(\mathfrak{c}_I|{\bf G}/\varphi(\partial {\bf e^*}))
:=\Prob_{{\bf e}\in {\bf \Omega}_I}[{\bf G}({\bf e})=\mathfrak{c}_I|
{\bf G}/\varphi(\partial {\bf e})
={\bf G}/\varphi(\partial {\bf e^*})
].\label{6O16a}
\end{eqnarray}
Denote by ${\bf reg}_{k,h,L}({\bf G})$ the minimum value 
of $\epsilon$ such that ${\bf G}$ is $(\epsilon,k,h,L)$-regular.
\end{df}
We will use the following new hypergraph regularity lemma \cite{I06}, which yields 
a shortest proof of Szemer\'edi's theorem on arithmetic progressions.
\begin{tha}[Regularity Lemma \cite{I06}]
\label{mt060721}
For any $r\ge k,h,\vec{b}=(b_i)_{i\in [k]},\epsilon>0,$ 
there exist integers 
$
\widetilde{m}_1,\cdots,\widetilde{m}_{k-1}
$ such that if 
 ${\bf G}$  is 
a $\vec{b}$-colored ($k$-bound $r$-partite hyper)graph on $
{\bf \Omega}$ then  for some 
integers $m_1,\cdots,m_{k-1}$ with 
$m_i\le \widetilde{m}_i, i\in [k-1],$ 
\begin{eqnarray*}
\Ex_{\vec{\varphi}\in\Phi(m_1,\cdots,m_{k-1})}
[{\bf reg}_{k,h}({\bf G}/\vec{\varphi})]
\le \epsilon.
\end{eqnarray*}
\end{tha}
The proof of the above in \cite{I06} essentially tells us the following.
\begin{tha}[Strong Form of Regularity Lemma \cite{I06}]
\label{6S19}
For any $r\ge k,h,\vec{b}=(b_i)_{i\in [k]},\epsilon>0,$ and 
for any function $L:\naturalset^{k-1}\to\naturalset$, 
there exist integers 
$\widetilde{m}_1,\cdots,\widetilde{m}_{k-1}
$ such that if 
 ${\bf G}$  is 
a $\vec{b}$-colored ($k$-bound $r$-partite hyper)graph on $
{\bf \Omega}$ then  for some 
integers $m_1,\cdots,m_{k-1}$ with 
$m_i\le \widetilde{m}_i, i\in [k-1],$ 
\begin{eqnarray*}
\Ex_{\vec{\varphi}\in\Phi(m_1,\cdots,m_{k-1})}
[{\bf reg}_{k,h,L(m_1,\cdots,m_{k-1})}({\bf G}/\vec{\varphi})]
\le \epsilon.
\end{eqnarray*}
\end{tha}
Theorem \ref{mt060721} is also used in \cite{I06lr} to show 
the hypergraph extension of the graph theorem by \cite{CRST}.
That is, the Ramsey number is linear (with respect to the order) 
for every bounded-degree hypergraph, which is also shown 
independently in Cooley et al. \cite{CFKO2} by a different way.
\mbox{}\bigskip\par
As I wrote in a final part of \cite{I06}, 
Property Testing and Regularization are essentially the same. 
They are all about random samplings, especially when considering 
constant-size (induced)subgraphs. 
If there exists a difference between the two, it is 
\begin{em}
whether the number of random vertex samplings 
is (PT) a fixed constant or (R) bounded by a constant but chosen randomly.
\end{em}
It may not be significant because a (non-canonical) 
property tester can invisualize some random number of 
vertex samples after choosing the vertices.\footnote{
Canonical property test  chooses (a fixed number of) 
vertices at random, but once the vertices are chosen, 
it outputs its answer deterministically.
Therefore, at first sight, 
canonical property tests may be weaker.
However as seen in \cite{AFKS}\cite[Th.2]{GT}, for any given non-canonical 
property test, there exists a canonical property test which is equivalent to it.
(Its derandomizing process is easy, since the sampling size of a non-canonical 
tester is a constant. The canonical tester repeats the samplings many (but constant) times.
 Then it computes the probability that the noncanonical tester accepts 
for each sampling. The canonical tester accepts iff 
the sum of the probabilities is at least $1/2.$
)}
\section{Lemmas and Their Proofs}
\begin{df}
Let ${\bf H}$ be a $k$-bound (colored 
$\mathfrak{r}$-partite hyper)graph on ${\bf \Omega}.$ 
Let $\delta: {\rm TC}({\bf H})\to [0,\infty)$ be a function.
Then for $I\in {\mathfrak{r}\choose [k]}, \alpha\in [0,\infty)$, 
we define a subset of  ${\rm TC}_I({\bf H})$ by 
\begin{eqnarray}
{\rm O}^\alpha_{\delta}
{\rm TC}_I({\bf H}):=\{
(\mathfrak{c}_J)_{J\subset I} |
{\bf d}_{\bf H}((\mathfrak{c}_J)_{J\subset I^*})\ge 
{\alpha^{1/3} \over |{\rm C}_{I^*}({\bf H})|} \mbox{ and } 
 \delta((\mathfrak{c}_J)_{J\subset I^*})\le 
{\alpha^{2/3} \over |{\rm C}_{I^*}({\bf H})|}
,
\forall I^*\subset I
\}.\label{..1227}
\end{eqnarray}
Write ${\rm \overline{O}}^\alpha_{\delta}
{\rm TC}_I({\bf H}):={\rm TC}_I({\bf H})\setminus 
{\rm O}^\alpha_{\delta}
{\rm TC}_I({\bf H})$.
(Here ${\rm O}$ means \lq ordinary\rq.)
We may drop the subscript $\delta$ if it is not necessary.
\par
Similarly we define ${\rm O}^\alpha_{\delta} \partial {\rm C}_I({\bf H})$ 
and ${\rm \overline{O}}^\alpha_{\delta} \partial {\rm C}_I({\bf H}).$
That is, $
{\rm O}^\alpha_{\delta}\partial 
{\rm C}_I({\bf H}):=\{
(\mathfrak{c}_J)_{J\subsetneq I} |
{\bf d}_{\bf H}((\mathfrak{c}_J)_{J\subset I^*})\ge 
{\alpha^{1/3} \over |{\rm C}_{I^*}({\bf H})|} \mbox{ and } 
 \delta((\mathfrak{c}_J)_{J\subset I^*})\le {\alpha^{2/3} \over |{\rm C}_{I^*}({\bf H})|}
,
\forall I^*\subsetneq I
\}.$
\end{df}
In the above notation, we easily see that if 
${\bf H}$ is $(\epsilon,k,1)$-regular then 
\begin{eqnarray}
\Prob_{{\bf e}\in {\bf \Omega}_I}[{\bf H}\ang{\bf e}\in 
{\rm \overline{O}}^{\epsilon}{\rm TC}_I({\bf H})
]
&\le& \sum_{J\subset I}
\Prob_{{\bf e}\in {\bf \Omega}_J}[
{\bf d}_{{\bf H}}({\bf H}\ang{\bf e})
<{\epsilon^{1/3}\over |{\rm C}_J({\bf H})|
}
\mbox{ or }
\delta({\bf H}\ang{\bf e})> {
\epsilon^{2/3}
\over |{\rm C}_J({\bf H})|
}
]
\nonumber\\
&\stackrel{(*)}{\le}
& \sum_{J\subset I} (\epsilon^{1/3}+\epsilon^{1/3})
\le 2^{|I|+1}\epsilon^{1/3} \label{6O16b}
\end{eqnarray}
where in the above (*) we used the fact that 
\begin{eqnarray*}
&&
\Prob_{{\bf e}\in {\bf \Omega}_J}\left[
\Prob_{{\bf e'}\in {\bf \Omega}_J}
[{\bf H}({\bf e'})={\bf H}({\bf e})
|{\bf e'}\stackrel{\partial {\bf H}}{\approx}{\bf e}
]
\le
{\epsilon^{1/3}
\over |{\rm C}_J({\bf H})|}
\right]
\nonumber
\\
&=&\sum_{\mathfrak{c}_J\in {\rm C}_J({\bf H})}
\Prob_{{\bf e}\in {\bf \Omega}_J}\left[
{\bf H}({\bf e})=\mathfrak{c}_J \mbox{ and }
\Prob_{{\bf e'}\in {\bf \Omega}_J}
[{\bf H}({\bf e'})
=\mathfrak{c}_J
|{\bf e'}\stackrel{\partial {\bf H}}{\approx}{\bf e}
]
\le
{\epsilon^{1/3}\over |{\rm C}_J({\bf H})|}
\right]\nonumber
\\
&\le &\sum_{\mathfrak{c}_J\in {\rm C}_J({\bf H})}
1\cdot
\Prob_{{\bf e}\in {\bf \Omega}_J}\left[
{\bf H}({\bf e})=\mathfrak{c}_J \left|
\Prob_{{\bf e'}\in {\bf \Omega}_J}
[{\bf H}({\bf e'})
=\mathfrak{c}_J
|{\bf e'}\stackrel{\partial {\bf H}}{\approx}{\bf e}
]
\le
{\epsilon^{1/3}\over |{\rm C}_J({\bf H})|}
\right.
\right]\nonumber
\\
&= &\sum_{\mathfrak{c}_J\in {\rm C}_J({\bf H})}
\Ex_{{\bf e}\in {\bf \Omega}_J}\left[
\Prob_{{\bf \tilde{e}}\in {\bf \Omega}_J}
[
{\bf H}({\bf \tilde{e}})=\mathfrak{c}_J 
|
{\bf \tilde{e}}\stackrel{\partial {\bf H}}{\approx}
{\bf e}
]
\left|
\Prob_{{\bf e'}\in {\bf \Omega}_J}
[{\bf H}({\bf e'})
=\mathfrak{c}_J
|{\bf e'}\stackrel{\partial {\bf H}}{\approx}{\bf e}
]
\le
{\epsilon^{1/3}\over |{\rm C}_J({\bf H})|}
\right.
\right]\nonumber
\\
&& \hspace{1.5in}\mbox{($\because$ the conditional part depends only on ${\bf H}(\partial {\bf e})$)}
\nonumber
\\
&\le  &\sum_{\mathfrak{c}_J\in {\rm C}_J({\bf H})}
\Ex_{{\bf e}\in {\bf \Omega}_J}\left[
{\epsilon^{1/3}\over |{\rm C}_J({\bf H})|}
\left|
\Prob_{{\bf e'}\in {\bf \Omega}_J}
[{\bf H}({\bf e'})
=\mathfrak{c}_J
|{\bf e'}\stackrel{\partial {\bf H}}{\approx}{\bf e}
]
\le
{\epsilon^{1/3}\over |{\rm C}_J({\bf H})|}
\right.
\right]
=\epsilon^{1/3}.
\end{eqnarray*}
\begin{df}[Color representative $\vartheta$]
\label{6S27}
Let ${\bf H}$ be a $k$-bound 
(colored $\mathfrak{r}$-partite hyper)graph on ${\bf \Omega}$.
Let $\vec{\psi}\in 
\Phi(\vec{m})$, where $\vec{m}\in \naturalset^{k-1},$ and $\epsilon,\epsilon_1>0.$
\par
For $\mathfrak{c}^*_I\in {\rm C}_I({\bf H}/\vec{\psi})$ 
with $I\in {\mathfrak{r}\choose [k]},$ 
we denote by ${\bf H}[\mathfrak{c}^*_I]$ the unique color $\mathfrak{c}_I
\in {\rm C}_I({\bf H})$ such that ${\bf H}/\vec{\psi}({\bf e})=\mathfrak{c}^*_I$ 
implies  ${\bf H}({\bf e})=\mathfrak{c}_I$.
Similarly we define ${\bf H}[\vec{\mathfrak{c}^*}]\in {\rm TC}_I({\bf H})$ 
for $\vec{\mathfrak{c}^*}\in {\rm TC}_I({\bf H}/\vec{\psi})$ 
and ${\bf H}[\mathfrak{c}^*]\in \partial {\rm C}_I({\bf H})$ 
 for ${\mathfrak{c}^*}\in \partial {\rm C}_I({\bf H}/\vec{\psi})$.
\par
Let $L_1,\cdots,L_{k}$ be positive integers.
Denote by ${\mathcal A}_I={\mathcal A}_I(L_1,\cdots,L_{|I|})$ 
the set of vectors $\vec{a}=(a_J)_{J\subset I}$ where 
$a_J\in [L_{|J|}].$ 
Write ${\mathcal A}_i:=\bigcup_{I\in {\mathfrak{r}\choose i}}{\mathcal A}_I$ 
and ${\mathcal A}:=\bigcup_{i\in [k]}{\mathcal A}_i$.
\par
$\bullet$ 
We inductively and probabilistically define colors 
$\mathfrak{d}(\vec{a})=\mathfrak{d}_I(\vec{a})\in {\rm C}_I({\bf H}/\vec{\psi})$ 
for all  $\vec{a}\in {\mathcal A}_I,I\in {\mathfrak{r}\choose [k]},$ 
by the following (i) and (ii).
\\
(i) Let $1\le s<k.$ Assume that 
$\mathfrak{d}_I(\vec{a})\in {\rm C}_I({\bf H}/\vec{\psi})$ is defined
 for any $I\in {\mathfrak{r}\choose [s-1]}$ and 
for any  $\vec{a}\in {\mathcal A}_I$ 
.
\\
(ii) Let $I\in {\mathfrak{r}\choose s}$ and  $\vec{a}\in {\mathcal A}_I.$ 
Pick an edge ${\bf e}\in {\bf \Omega}_I$ randomly so that 
${\bf H}/\vec{\psi}(\partial {\bf e})=
(\mathfrak{d}_{J}(\vec{a}|_J))_{J\subsetneq I}$ 
where $\vec{a}|_J:=(a_{J'})_{J'\subset J}$.
Let $\mathfrak{d}_I(\vec{a}):={\bf H}/\vec{\psi}({\bf e}).$
\par
Note that for the entire process we pick a random edge 
exactly $|{\mathcal A}|=rL_1+{r\choose 2}L_1^2L_2+\cdots+{r\choose k}
\prod_{i\in [k]}L_i^{k\choose i}$ times.
\par
Write
$$
\vec{\mathfrak{d}}(\vec{a})=
\vec{\mathfrak{d}}_I(\vec{a})
:=(\mathfrak{d}_J(\vec{a}|_J))_{J\subset I}\in {\rm TC}_I({\bf H}/\vec{\psi})
\mbox{ and }
\vec{\mathfrak{d}}(\partial \vec{a})
=
\vec{\mathfrak{d}}_I(\partial \vec{a})
:=(\mathfrak{d}_J(\vec{a}|_J))_{J\subsetneq I}\in
 \partial {\rm C}_I({\bf H}/\vec{\psi})
$$
where $\partial \vec{a}:=
(a_{J'})_{J'\subsetneq I}.$
\par
$\bullet$ Assume that $\mathfrak{d}$ is fixed. 
Then we will inductively and probabilistically 
define a map $\theta_I:{\rm TC}_I({\bf H})\to [0,L_{|I|}]$ 
for all $I\in {\mathfrak{r}\choose [k]},$ 
by the following (i') and (ii').
\\
(i') Let $1\le s\le k.$ Assume that 
$\theta_I(\vec{\mathfrak{c}})\in [0,L_{|I|}]$ is defined
 for any $I\in {\mathfrak{r}\choose [s-1]}$ and 
for any  $\vec{\mathfrak{c}}\in {\rm TC}_I({\bf H})$.
\\
(ii') Let $I\in {\mathfrak{r}\choose s}$ and  $\vec{\mathfrak{c}}
=(\mathfrak{c}_J)_{J\subset I}\in {\rm TC}_I({\bf H}).$ 
Let 
$a_J:=\theta_J((\mathfrak{c}_{J'})_{J'\subset J})$ for $J\subsetneq I$.
If $a_J=0$ for some $J\subsetneq I$ then we define 
$\theta_I(\vec{\mathfrak{c}}):=0$. 
Suppose that $a_J\in [L_{|J|}]$ for all $J\subsetneq I.$
Let $$L_I^*:=
\{a_I\in [L_{|I|}] |
{\bf H}[\mathfrak{d}((a_J)_{J\subset I})]=\mathfrak{c}_I
\}.
$$ 
If $L_I^*\not=\emptyset$ then 
we define $\theta_I(\vec{\mathfrak{c}}):=a_I$ for 
an $a_I\in L_I^*$ chosen uniformly at random.
If $L_I^*=\emptyset$ then we define $\theta_I(\vec{\mathfrak{c}}):=0.$ 
\par
Write
$$
\vec{\theta}(\vec{\mathfrak{c}})=\vec{\theta}_I(\vec{\mathfrak{c}})
:=(\theta_J(\vec{\mathfrak{c}}|_J))_{J\subset I}
\mbox{ and }
\vec{\theta}(\partial \vec{\mathfrak{c}})
=\vec{\theta}_I(\partial \vec{\mathfrak{c}})
:=(\theta_J(\vec{\mathfrak{c}}|_J))_{J\subsetneq I}
$$
where $\vec{\mathfrak{c}}|_J:=(\mathfrak{c}_{J'})_{J'\subset J}$ and 
$\partial \vec{\mathfrak{c}}:=
(\mathfrak{c}_{J'})_{J'\subsetneq I}.$
\par
$\bullet$ When $\vec{\theta}_I(\vec{\mathfrak{c}})\in {\mathcal A}_I$ 
or $\vec{\theta}_I(\partial \vec{\mathfrak{c}})\in \partial {\mathcal A}_I
:=\{\vec{a}=(a_J)_{J\subsetneq I}| a_J\in [L_{|J|}]\}$ 
(i.e. the case when it does not contain any zero), 
we write 
$$
{\vartheta}_I(\vec{\mathfrak{c}}):=
{\mathfrak{d}}(
\vec{\theta}_I(\vec{\mathfrak{c}}))\in {\rm C}_I({\bf H}/\vec{\psi}),
\vec{\vartheta}_I(\vec{\mathfrak{c}}):=
\vec{\mathfrak{d}}(
\vec{\theta}_I(\vec{\mathfrak{c}}))\in {\rm TC}_I({\bf H}/\vec{\psi})
\mbox{ or }
\vec{\vartheta}_I(\partial \vec{\mathfrak{c}}):=
\vec{\mathfrak{d}}(
\vec{\theta}_I(\partial \vec{\mathfrak{c}}))\in \partial {\rm C}_I({\bf H}/\vec{\psi})
.$$
Otherwise, write ${\vartheta}_I(\vec{\mathfrak{c}}):=\mathfrak{0}$, 
$\vec{\vartheta}_I(\vec{\mathfrak{c}}):=\mathfrak{0}$ and 
$\vec{\vartheta}_I(\partial \vec{\mathfrak{c}}):=\mathfrak{0}$, 
where $\mathfrak{0}$ is a fixed symbol which does not belong to any color class.
\end{df}
In the proofs, we will 
 write ${\bf d}_{\bf G}^{(\delta)}(\vec{\mathfrak{c}})
={\bf d}_{\bf G}(\vec{\mathfrak{c}})\dot{\pm}\delta(\vec{\mathfrak{c}})$ for 
$\vec{\mathfrak{c}}\in {\rm TC}({\bf G}).$
\begin{lemma}[All representatives are very regular]
\label{z0813}
There exist a positive-valued function $\epsilon^{(\ref{z0813})}_1(\cdots)$ 
such that the following proposition holds.\par
Let $r\ge k$ be positive integers and let $\vec{L}=(L_i)_{i\in [k]}$ be 
a sequence of positive integers.
Let 
\begin{eqnarray*}
0<\epsilon_1\le \epsilon^{(\ref{z0813})}_1(r,k,\vec{L})
\end{eqnarray*}
and 
 ${\bf H}$  a $k$-bound 
(colored 
$\mathfrak{r}$-partite hyper)graph on ${\bf \Omega}.$
Suppose 
that ${\bf H}/\vec{\psi}$ is $(\epsilon_1,k,1)$-regular 
for some $\vec{\psi}=(\psi_i)_{i\in [k-1]}\in 
\Phi(m_1,\cdots,m_{k-1}),$ where $m_1,\cdots,m_{k-1}$ are positive integers.
Then the $\vartheta$ probabilistically defined 
 in Definition \ref{6S27} satisfy the following inequality:\\
\begin{eqnarray*}
\sum_{I\in {\mathfrak{r}\choose [k]}}
\sum_{\vec{\mathfrak{c}}\in 
{\rm TC}_I({\bf H})}
\Prob_{\vartheta}\left[
\vec{\vartheta}_I(\vec{\mathfrak{c}})
\in 
{\rm \overline{O}}^{\epsilon_1}{\rm TC}_I({\bf H}/
\vec{\psi})
\right]
<  0.01,
\end{eqnarray*}
where 
$\Prob_{\vartheta}$ denotes the probability in the probability space 
generated by the (two-step) random process in the definition of $\vartheta$.
\end{lemma}
\Proof
By the regularity of ${\bf H}/\vec{\psi}$, 
we see that 
\begin{eqnarray}
&&
\sum_{I\in {\mathfrak{r}\choose [k]}}
\sum_{\vec{\mathfrak{c}}\in 
{\rm TC}_I({\bf H})}
\Prob_{\vartheta}\left[
\vec{\vartheta}(\vec{\mathfrak{c}})
\in 
{\rm \overline{O}}^{\epsilon_1}
{\rm TC}_I({\bf H}/
\vec{\psi})
\right]\nonumber
\\
&\le &
\sum_{I\in {\mathfrak{r}\choose [k]}}
\sum_{\vec{a}\in {\mathcal A}_I}
\Prob_{\mathfrak{d}}[
\vec{\mathfrak{d}}_I(\vec{a})\in 
{\rm \overline{O}}^{\epsilon_1}{\rm TC}_I({\bf H}/\vec{\psi})
]
\quad (\because 
{\vec{\theta}}_I(\vec{\mathfrak{c}}),
{\vec{\theta}}_I(\vec{\mathfrak{c}'})\in {\mathcal A}_I, 
\vec{\mathfrak{c}}\not=\vec{\mathfrak{c}'}
\mbox{ implies } 
{\vec{\theta}}_I(\vec{\mathfrak{c}})\not=
{\vec{\theta}}_I(\vec{\mathfrak{c}'})
)
\nonumber\\
&=&
\sum_{I\in {\mathfrak{r}\choose [k]}}
\sum_{\vec{a}\in {\mathcal A}_I}
\brkt{
1-\sum_{\vec{\mathfrak{c}}^*\in 
{\rm O}^{\epsilon_1}{\rm TC}_I({\bf H}/\vec{\psi})}
\Prob_{\mathfrak{d}}[
\vec{\mathfrak{d}}_I(\vec{a})=\vec{\mathfrak{c}}^*
]
}
\nonumber\\
&=&
\sum_{I\in {\mathfrak{r}\choose [k]}}
|{\mathcal A}_I|
\brkt{
1-\sum_{\vec{\mathfrak{c}}^*\in 
{\rm O}^{\epsilon_1}{\rm TC}_I({\bf H}/\vec{\psi})}
\prod_{J\subset I}
\Prob_{{\bf e}\in {\bf \Omega}_J}[{\bf H}/\vec{\psi}({\bf e})=\mathfrak{c}^*_J
|
{\bf H}/\vec{\psi}(\partial {\bf e})=(\mathfrak{c}^*_{J'})_{J'\subsetneq J}
]
}
\nonumber\\
&=&
\sum_{I\in {\mathfrak{r}\choose [k]}}
|{\mathcal A}_I|
\brkt{
1-\sum_{\vec{\mathfrak{c}}^*\in 
{\rm O}^{\epsilon_1}{\rm TC}_I({\bf H}/\vec{\psi})}
{
\Prob_{{\bf e}\in {\bf \Omega}_I}[{\bf H}/\vec{\psi}\ang{\bf e}
=\vec{\mathfrak{c}}^*]
\over 
\prod_{J\subset I}
{\bf d}^{(\delta)}_{{\bf H}/\vec{\psi}}((\mathfrak{c}^*_{J'})_{J'\subset J})}
\cdot 
\prod_{J\subset I}
{\bf d}_{{\bf H}/\vec{\psi}}((\mathfrak{c}^*_{J'})_{J'\subset J})
}
\nonumber\\
&& 
\quad (\because \mbox{$\vec{\mathfrak{c}}^*$ can be considered as 
a simplicial-complex in } {\mathcal S}_{r,|I|,1,{\bf H}/\vec{\psi}}.
\mbox{ Use regularity (i) of }{\bf H}/\vec{\psi}.
)
\nonumber\\
&\stackrel{(\ref{..1227})}{=}&
\sum_{I\in {\mathfrak{r}\choose [k]}}\brkt{
\prod_{j\in [|I|]}L_j^{|I|\choose j}
}
\brkt{
1-{1\over \brkt{1\dot{\pm}\epsilon_1^{1/3}}^{2^{|I|}}}
\Prob_{{\bf e}\in {\bf \Omega}_I}[{\bf H}/\vec{\psi}\ang{\bf e}
\in {\rm O}^{\epsilon_1}{\rm TC}_I({\bf H}/\vec{\psi})
]
}
\nonumber\\
&\stackrel{(\ref{6O16b})}{\le} &
\sum_{I\in {\mathfrak{r}\choose [k]}}
\brkt{
\prod_{j\in [|I|]}L_j^{|I|\choose j}
}
\brkt{
1-{1-2^{|I|+1} \epsilon_1^{1/3}\over \brkt{1+\epsilon_1^{1/3}}^{2^{|I|}}}
}
\nonumber\\
&< &0.01\nonumber
\end{eqnarray}
where the last inequality follows from the assumption that 
$\epsilon_1>0$ is small enough with respect to $r,k,L_1,\cdots,L_{k}.$ 
\lemmaqed
\begin{lemma}[Most representatives are ordinary]
\label{6O12d}
There exist positive-valued 
functions $\epsilon^{(\ref{6O12d})}(\cdot), 
\epsilon_1^{(\ref{6O12d})}(\epsilon)=c\epsilon^2,$
where $c>0$ is a small absolute constant, and 
$L^{(\ref{6O12d})}=L(\cdot,\cdot)$ such that 
the following proposition holds.
\par
Let $r\ge k$ be positive integers. Let $0<\epsilon\le \epsilon^{(\ref{6O12d})}(k)$ and 
$0<\epsilon_1\le \epsilon_1^{(\ref{6O12d})}(\epsilon)$.
Let $\vec{b'}=(b'_i)_{i\in [k]}$ and 
 $(L_i)_{i\in [k-1]}$ be 
sequences of positive integers with $L_i\ge L(\epsilon,b'_i)$ for all $i<k.$
Let ${\bf H}$ be a $k$-bound $\vec{b'}$-colored 
($\mathfrak{r}$-partite hyper)graph on ${\bf \Omega}$.
For some integers $m_1,\cdots,m_{k-1}$, we 
suppose that ${\bf H}$ is 
$(\epsilon,k,1,m_1+\cdots+m_{k-1})$-regular and 
that 
\begin{eqnarray}
\Ex_{\vec{\psi}\in \Phi(m_1,\cdots,m_{k-1})}
[{\bf reg}_{k,1}({\bf H}/\vec{\psi})]\le \epsilon_1^2.
\label{6O12}
\end{eqnarray}
Then the $\vartheta$ probabilistically defined 
in Definition \ref{6S27} satisfies the following :
\begin{eqnarray*}
\Ex_{\vec{\psi}\in \Phi(m_1,\cdots,m_{k-1})}
\Ex_{\mathfrak{\vartheta}}\Prob_{{\bf e}\in {\bf \Omega}_I}[
\vec{\vartheta}({\bf H}(\partial {\bf e}))
\mbox{ is $(\epsilon_1,\epsilon^{2/3},\epsilon^{1/3})$-ordinary }
]\ge 1-2^{|I|}\epsilon^{1/3}/c
\quad \forall I\in {\mathfrak{r}\choose [k]},
\end{eqnarray*}
where we call $\mathfrak{c}^*=(\mathfrak{c}^*_J)_{J\subsetneq I}\in 
\partial {\rm C}_I({\bf H}/\vec{\psi})$  {\bf $(\epsilon_1,\gamma,\alpha)$-ordinary}
 if and only if \\
(i) $\mathfrak{c}^*
\in {\rm O}^{\epsilon_1}
\partial {\rm C}_I({\bf H}/\vec{\psi})$, \\
(ii)
for all $J\subset I$
\begin{eqnarray}
\sum_{\mathfrak{c}_J\in {\rm C}_J({\bf H})}
\left|
{\bf d}_{{\bf H}/\vec{\psi}}(\mathfrak{c}_J|(\mathfrak{c}^*_{J'})_{J'\subsetneq J})
-
{\bf d}_{{\bf H}}(\mathfrak{c}_J|
({\bf H}[\mathfrak{c}^*_{J'}])_{J'\subsetneq J}
)
\right|^2
\le \brkt{\gamma\over |{\rm C}_J({\bf H})|}^2, 
\mbox{ and }
\label{..1227a}
\end{eqnarray}
\\
(iii) if $\mathfrak{c}_I\in {\rm C}_I({\bf H})$ and 
${\bf d}_{{\bf H}}(\mathfrak{c}_I|({\bf H}[\mathfrak{c}^*_{J'}])_{J'\subsetneq I})
\ge 
{\alpha
\over |{\rm C}_I({\bf H})|
}$ then $\vec{\vartheta}((\mathfrak{c}_J)_{J\subset I})\not=\mathfrak{0}.$
\end{lemma}
In the above, we mean 
$
{\bf d}_{{\bf H}/\vec{\psi}}(\mathfrak{c}_J|(\mathfrak{c}^*_{J'})_{J'\subsetneq J})
:=\Prob_{{\bf e}\in {\bf \Omega}_I}[{\bf G}({\bf e})=\mathfrak{c}_I|
{\bf G}/\vec{\psi}(\partial {\bf e})=(\mathfrak{c}^*_{J'})_{J'\subsetneq J}
]
$ as in (\ref{6O16a}).
\\
\Proof 
In the below, we write ${\bf H}^*:={\bf H}/\vec{\psi}.$
Let $\gamma,\rho>0$, which will be defined later at (\ref{6O11b}).
Write 
\begin{eqnarray*}
{\rm O}^*\partial {\rm C}_I({\bf H}^*)
:=\left\{\left.
\mathfrak{c}\in 
{\rm O}^{\epsilon_1}
\partial {\rm C}_I({\bf H}^*)
\right| \mathfrak{c} \mbox{ is $(\epsilon_1,\gamma,\infty)$-ordinary}
\right\}.
\end{eqnarray*}
We say that $\mathfrak{c}\in \partial {\rm C}_I({\bf H})$ is 
{\bf $(\epsilon_1,\gamma;\rho)$-ordinary} if and only if 
$\Prob_{{\bf e^*}\in {\bf \Omega}_I}[{\bf H}^*(\partial {\bf e^*})
\not\in {\rm O}^*\partial {\rm C}_I({\bf H}^*)
|{\bf H}(\partial {\bf e^*})=\mathfrak{c}
]\le \rho.$
Write 
\begin{eqnarray}
{\rm O}^\circ \partial {\rm C}_I({\bf H})
:=\{\mathfrak{c}\in 
{\rm O}^{\epsilon}\partial {\rm C}_I({\bf H})
| \mathfrak{c} \mbox{ is $(\epsilon_1,\gamma;\rho)$-ordinary}
\}.\label{..1227b}
\end{eqnarray}
Since 
$(\epsilon,k,1,m_1+\cdots+m_{k-1})$-regularity of ${\bf H}$ 
yields that 
$$
\Prob_{\vec{\psi}\in\Phi(m_1+\cdots+m_{k-1}),
{\bf e^*}\in {\bf \Omega}_J}[\sum_{\mathfrak{c}_J\in {\rm C}_J({\bf H})
}
\brkt{
{\bf d}_{{\bf H}/\varphi}(\mathfrak{c}_J|{\bf H}/
\varphi(\partial {\bf e^*}))
-
{\bf d}_{{\bf H}}(\mathfrak{c}_J|{\bf H}(\partial {\bf e^*}))
}^2
]\le \brkt{\epsilon\over |{\rm C}_J({\bf H})|}^2
$$ for all $J\in {\mathfrak{r}\choose [k]}$, i.e. 
(by the definition of regularization)
$$
\Prob_{\vec{\psi}\in\Phi(\vec{m}),
{\bf e^*}\in {\bf \Omega}_I}[\sum_{\mathfrak{c}_J\in {\rm C}_J({\bf H})
}
\brkt{
{\bf d}_{{\bf H}/\vec{\psi}}(\mathfrak{c}_J
|{\bf H}/\vec{\psi}(\partial ({\bf e^*}|_J)))
-
{\bf d}_{{\bf H}}(\mathfrak{c}_J|{\bf H}(\partial ({\bf e^*}|_J)))
}^2
]\le \brkt{\epsilon\over |{\rm C}_J({\bf H})|}^2
$$ for all $J\subset I$
,
 it is easy to see that 
\begin{eqnarray}
\underset{\vec{\psi}\in\Phi(\vec{m}),
{\bf e^*}\in {\bf \Omega}_I}{\Prob}
[
%
\sum_{\mathfrak{c}_I\in {\rm C}_J({\bf H})
}
\brkt{
{\bf d}_{{\bf H}/\vec{\psi}}(\mathfrak{c}_J|
{\bf H}/\vec{\psi}(\partial ({\bf e^*}|_J)))
-
{\bf d}_{{\bf H}}(\mathfrak{c}_J|{\bf H}(\partial ({\bf e^*}|_J)))
}^2
> \brkt{\gamma\over |{\rm C}_J({\bf H})|}^2
\exists J\subset I
%
]
\le 2^{|I|}{\epsilon^2\over \gamma^2},
\label{.320}
\end{eqnarray}
which yields that 
\begin{eqnarray}
&&
\Ex_{\vec{\psi}\in\Phi(\vec{m})}
\Prob_{{\bf e^*}\in {\bf \Omega}_I}
[{\bf H}^*(\partial {\bf e^*})\not\in 
{\rm O}^*\partial {\rm C}_I({\bf H}^*)]
\nonumber\\
&\le&
(\ref{.320})
+\Prob_{\vec{\psi}\in\Phi(\vec{m}),
{\bf e^*}\in {\bf \Omega}_I}
[{\bf H}^*(\partial {\bf e^*})\in 
{\rm \overline{O}}^{\epsilon_1}
\partial {\rm C}_I({\bf H}^*)]
\nonumber\\
&\stackrel{(\ref{6O16b})}{\le}&
2^{|I|}\epsilon^2/\gamma^2
+
\sum_{J\subsetneq I}
(\epsilon_1^{1/3}+\epsilon_1^{1/3})+
\Prob_{\vec{\psi}}[{\bf reg}_{k,1}({\bf H}^*)> \epsilon_1]
j\nonumber\\
&\stackrel{(\ref{6O12})}{<}
& 2^{|I|}({\epsilon^2\over \gamma^2}+3\epsilon_1^{1/3}).\label{6O16c}
\end{eqnarray}
Thus we see that
\begin{eqnarray}
&&
\Ex_{\vec{\psi}\in\Phi(\vec{m})}
\Prob_{{\bf e}\in {\bf \Omega}_I}
[{\bf H}(\partial {\bf e})\not\in 
{\rm O}^\circ \partial {\rm C}_I({\bf H})
]
\nonumber\\
&\le& 
\Prob_{\vec{\psi}\in\Phi(\vec{m}),
{\bf e}\in {\bf \Omega}_I}
[{\bf H}(\partial {\bf e}) \mbox{ is not $(\epsilon_1,\gamma;\rho)$-ordinary}]
+\Prob_{{\bf e}\in {\bf \Omega}_I}
[{\bf H}(\partial {\bf e})\in {\rm \overline{O}}^{\epsilon}
\partial {\rm C}_I({\bf H})]
\nonumber\\
&\le& 2^{|I|}
\brkt{
{\epsilon^2/\gamma^2+3\epsilon_1^{1/3}
\over  \rho}+2\epsilon^{1/3}}.
 \hspace{1cm}(\because (\ref{6O16c}),(\ref{6O16b}))
\label{6O11a}
\end{eqnarray}
Therefore 
if $\mathfrak{c}
\in 
{\rm O}^{\epsilon}\partial {\rm C}_I({\bf H})$ and 
 $\mathfrak{c}^*
\in 
{\rm O}^*\partial {\rm C}_I({\bf H}^*)$
then 
\begin{eqnarray}
&&
\Prob_{\mathfrak{d},\theta}[\vec{\vartheta}(\mathfrak{c})=\mathfrak{c}^*
]
\nonumber\\
&= &
\prod_{J\subsetneq I}
{
\Prob_{{\bf e}\in {\bf \Omega}_J}
[{\bf H}^*({\bf e})=\mathfrak{c}^*_J|
{\bf H}^*(\partial {\bf e})=(\mathfrak{c}^*_{J'})_{J'\subsetneq J}]
\over 
\Prob_{{\bf e}\in {\bf \Omega}_J}
[{\bf H}({\bf e})=\mathfrak{c}_J|
{\bf H}^*(\partial {\bf e})=(\mathfrak{c}^*_{J'})_{J'\subsetneq J}]
}
\brkt{1-
(1-\Prob_{{\bf e}\in {\bf \Omega}_J}
[{\bf H}({\bf e})=\mathfrak{c}_J|
{\bf H}^*(\partial {\bf e})=(\mathfrak{c}^*_{J'})_{J'\subsetneq J}]
)^{L_{|J|}}}
\nonumber\\
&\stackrel{(\ref{..1227a})}{=}&
\prod_{J\subsetneq I}
{
{\bf d}_{{\bf H}^*}((\mathfrak{c}^*_{J'})_{J'\subset J})
\over 
{\bf d}_{{\bf H}}((\mathfrak{c}_{J'})_{J'\subset J})
\dot{\pm}\gamma/|{\rm C}_J({\bf H})|
}
\brkt{1-
(1-
{\bf d}_{{\bf H}}((\mathfrak{c}_{J'})_{J'\subset J})
\dot{\mp}\gamma/|{\rm C}_J({\bf H})|
)^{L_{|J|}}
}
\nonumber\\
&\ge &{
\Prob_{{\bf e}\in {\bf \Omega}_I}[{\bf H}^*(\partial {\bf e})=
\mathfrak{c}^*]
\over 
\Prob_{{\bf e}\in {\bf \Omega}_I}[{\bf H}(\partial {\bf e})=
\mathfrak{c}]
}
\prod_{J\subsetneq I}
{
{\bf d}^{(\delta)}_{{\bf H}}((\mathfrak{c}_{J'})_{J'\subset J})
\over 
{\bf d}^{(\delta)}_{{\bf H}^*}((\mathfrak{c}^*_{J'})_{J'\subset J})
}
{
{\bf d}_{{\bf H}^*}((\mathfrak{c}^*_{J'})_{J'\subset J})
\over 
{\bf d}_{{\bf H}}((\mathfrak{c}_{J'})_{J'\subset J})
\dot{\pm}\gamma/|{\rm C}_J({\bf H})|
}
\brkt{1-
\brkt{1-
{\epsilon^{1/3}-\gamma
\over |{\rm C}_J({\bf H})|
}
}^{L_{|J|}}
}
\nonumber\\ &&
(\mbox{Use regularities where }
\mathfrak{c},\mathfrak{c}^* \mbox{ are considered complexes in }
{\mathcal S}_{r,|I|-1,1,{\bf H}} \mbox{ and in } {\mathcal S}_{r,|I|-1,1,{\bf H^*}}.
)
\nonumber\\
&\ge &{
\Prob_{{\bf e}\in {\bf \Omega}_I}[{\bf H}^*(\partial {\bf e})=
\mathfrak{c}^*
|
{\bf H}(\partial {\bf e})=
\mathfrak{c}]
}
\brkt{
{
1\dot{\pm}\epsilon^{1/3}
\over 
(1\dot{\pm}\epsilon_1^{1/3})
}
{1\over (1\dot{\pm}{\gamma\over \epsilon^{1/3}})}
\cdot (1-\epsilon)
}^{2^{|I|}}.
\quad (\because L_{|J|}\ge L(\epsilon,b'_{|J|}).)
\label{6O11}
\end{eqnarray}
Hence it follows that 
\begin{eqnarray}
&&\Ex_{\vec{\psi}}
\Ex_{{\vartheta}}
\Prob_{{\bf e}\in {\bf \Omega}_I}[
\vec{\vartheta}({\bf H}(\partial {\bf e}))
\mbox{ is  $(\epsilon_1,\gamma,\infty)$-ordinary}
]
\nonumber\\
&\ge &\Ex_{\vec{\psi}}
\sum_{\mathfrak{c}\in {\rm O}^\circ \partial {\rm C}_I({\bf H})}
\sum_{\mathfrak{c}^*\in {\rm O}^*\partial {\rm C}_I({\bf H^*})}
\Prob_{{\bf e}\in {\bf \Omega}_I}[{\bf H}(\partial {\bf e})=\mathfrak{c}
]
\Prob_{{\vartheta}}
[
\vec{\vartheta}(\mathfrak{c})=\mathfrak{c}^*
]
\nonumber\\
&
\stackrel{(\ref{6O11})}{\ge} &
\brkt{
{
1-\epsilon^{1/3}
\over 
(1+\epsilon_1^{1/3})
}
{1-\epsilon\over (1+{\gamma\over \epsilon^{1/3}})}
}^{2^{|I|}}\Ex_{\vec{\psi}}
\sum_{\mathfrak{c}\in {\rm O}^\circ \partial {\rm C}_I({\bf H})}
\Prob_{{\bf e}\in {\bf \Omega}_I}[{\bf H}(\partial {\bf e})=\mathfrak{c}
]
\sum_{\mathfrak{c}^*}
\Prob_{{\bf e}\in {\bf \Omega}_I}[{\bf H}^*(\partial {\bf e})=
\mathfrak{c}^*|{\bf H}(\partial {\bf e})=\mathfrak{c}
]
\nonumber\\
&\stackrel{(\ref{..1227b})}{\ge} &
\brkt{
{
1-\epsilon^{1/3}
\over 
(1+\epsilon_1^{1/3})
}
{1-\epsilon\over (1+{\gamma\over \epsilon^{1/3}})}
}^{2^{|I|}}
(1-\rho)\Ex_{\vec{\psi}}
\sum_{\mathfrak{c}}
\Prob_{{\bf e}\in {\bf \Omega}_I}[{\bf H}(\partial {\bf e})=\mathfrak{c}
]
\nonumber\\
&
\stackrel{(\ref{6O11a})}{\ge} &
\brkt{
{
1-\epsilon^{1/3}
\over 
(1+\epsilon_1^{1/3})
}
{1-\epsilon\over (1+{\gamma\over \epsilon^{1/3}})}
}^{2^{|I|}}
(1-\rho)
\brkt{1-2^{|I|}
\brkt{{\epsilon^2/\gamma^2+3\epsilon_1^{1/3}
\over  \rho}+2\epsilon^{1/3}}
}
\nonumber\\
&{\ge} &
1-
1.1 \brkt{
 2^{|I|}(3\epsilon^{1/3}+{\gamma\over \epsilon^{1/3}}
+
{\epsilon^2/\gamma^2+3\epsilon_1^{1/3}
\over  \rho}
)
+\rho}
\nonumber\\
&{\ge} &
1-
1.1 
(5.1\cdot 2^{|I|}+1)
\epsilon^{1/3}
\hspace{1cm} (\mbox{where }\gamma:=\epsilon^{2/3},\rho:=\epsilon^{1/3},\mbox{ and }
\epsilon_1=o(\epsilon^2))
.\label{6O11b}
\end{eqnarray}
Finally, we have that 
\begin{eqnarray}
&&\Ex_{\vec{\psi}}
\Ex_{{\vartheta}}
\Prob_{{\bf e}\in {\bf \Omega}_I}[
\vec{\vartheta}({\bf H}(\partial {\bf e}))
\mbox{ is $(\epsilon_1,\epsilon^{2/3},\epsilon^{1/3})$-ordinary}
|\vec{\vartheta}({\bf H}(\partial {\bf e}))
\mbox{ is  $(\epsilon_1,\epsilon^{2/3},\infty)$-ordinary}]
\nonumber\\
&\ge &1-
\sum_{\mathfrak{c}_I\in 
{\rm C}_I({\bf H})}
\brkt{1-{\epsilon^{1/3}
\over |{\rm C}_I({\bf H})|}}^{L_{|I|}}
\nonumber\\
&\ge& 1-\epsilon^{1/3},\label{6O15}
\end{eqnarray}
when $L_{|I|}\ge L(\epsilon,b'_{|I|})$.
Combining (\ref{6O11b}) and (\ref{6O15}) completes the proof.
\lemmaqed
\begin{df}[Abbreviation]
Let ${\bf G}$ be a $k$-bound $\vec{b}=(b_i)_{i\in [k]}$-colored hypergraph.
Write $c_i({\bf G}):=\max_{I\in {{\mathfrak r}\choose i}}|{\rm C}_I({\bf G})|$ 
for $i\in [k]$.
For an integer $m$, we write 
$\vec{B}(\vec{b},m):=({B}_i(\vec{b},m))_{i\in [k]}$ 
where $
B_i(\vec{b},m):=
\prod_{j\in [0,k-i]}b_{i+j}^{{r\choose j}m^j}.
$ Note that 
\begin{eqnarray}
c_i({\bf G}/\varphi)\le
B_i(\vec{b},m)\quad \forall i\in [k]
\forall \varphi\in\Phi(m).\label{6O16d}
\end{eqnarray}
\end{df}
\begin{lemma}[Main Lemma]\label{6S14f}
There exists a positive-valued function $\epsilon_1^{(\ref{6S14f})}(\cdots)$
 such that 
the following proposition holds.\par
Let $r\ge k,\vec{b}=(b_i)_{i\in [k]}$ be positive integers and 
$0<\epsilon\le \epsilon^{(\ref{6O12d})}(k)$, where 
$\epsilon^{(\ref{6O12d})}(\cdot)$ is the function of Lemma \ref{6O12d}.
Let $\epsilon_1^\circ:\naturalset^{k-1}\to (0,1]$ be a function
such that 
$$0<\epsilon_1^\circ(l_1,\cdots,l_{k-1})\le 
\epsilon_1^{(\ref{6S14f})}(r,k,\vec{b},\epsilon;l_1,\cdots,l_{k-1})
$$ for all integers $l_1,\cdots,l_{k-1}.$
Let $h^\circ:\naturalset^{k-1}\to\naturalset$ be a function.
\par
Then there exist an integer $\widetilde{\ell}$ and a 
function $\widetilde{m}(\cdots)$ 
such that 
if ${\bf G}$ is a $\vec{b}$-colored ($k$-bound $\mathfrak{r}$-partite hyper)graph 
(on ${\bf \Omega}$) then there exist integers 
$\ell_1,\cdots,\ell_{k-1}\in [\widetilde{\ell}]$ 
and integers $m_1,\cdots,m_{k-1}\in [\widetilde{m}(\ell_1,\cdots,\ell_{k-1})]$ 
which satisfy 
 the following, where $\epsilon_1:=\epsilon_1^\circ(
\ell_1,\cdots,\ell_{k-1})$.
\par
There exist $\vec{\varphi}\in \Phi(\ell_1,\cdots,\ell_{k-1})$ and 
$\vec{\psi}\in \Phi(m_1,\cdots,m_{k-1})$ 
such that  
(${\bf H}:={\bf G}/\vec{\varphi}$ is 
$(\epsilon,k,1,m_1+\cdots+m_{k-1})$-regular and that)
\begin{eqnarray*}
\mbox{
 ${\bf H}/\vec{\psi}$ is 
$(\epsilon_1,k,h^\circ(c_1({\bf H}),\cdots,c_{k-1}({\bf H})
))$-regular,}
\end{eqnarray*} 
and furthermore that the map 
$\vec{\vartheta}(\cdot)$ defined in Definition \ref{6S27} 
for the 
${\bf H}$ and $\vec{\psi}$ with some integers $(L_i)_{i\in [k]}$ 
satisfies 
all of the following properties for all $I\in {\mathfrak{r}\choose [k]}$,
simultaneously, 
with probability at least 0.9.\\ 
(i) 
$
\Prob_{{\bf e}\in {\bf \Omega}_I}[
\vec{\vartheta}({\bf H}(\partial {\bf e}))
\mbox{ is $(\epsilon_1,\epsilon^{2/3},\epsilon^{1/3})$-ordinary }
]\ge 1-O_{r,k}(\epsilon^{1/3}).$\\
(ii) 
If $\vec{\mathfrak{c}}\in {\rm TC}_I({\bf H})$ 
and $\vec{\vartheta}(\vec{\mathfrak{c}})\not=\mathfrak{0}$ then 
$
\vec{\vartheta}(\vec{\mathfrak{c}})
\in {\rm O}^{\epsilon_1}
{\rm TC}_I({\bf H}/\vec{\psi}).$
\\
(iii) 
If 
${\mathfrak{c}}=(\mathfrak{c}_J)_{J\subsetneq I}
\in \partial {\rm C}_I({\bf H})$ and 
$\vec{\vartheta}(\mathfrak{c})\not=\mathfrak{0}$ then 
there exists a color $\mathfrak{c}_I\in {\rm C}_I({\bf H})$ such that 
$\vec{\vartheta}((\mathfrak{c}_J)_{J\subset I})
\in {\rm O}^{\epsilon_1}
{\rm TC}_I({\bf H}/\vec{\psi}).$
\end{lemma}
\Proof Fix $r\ge k,\vec{b},\epsilon, \epsilon_1^\circ,
h^\circ$ and ${\bf G}$ as in the lemma.
Without loss of generality, $h^\circ$ is increasing.
The upper bound function $\epsilon_1^{(\ref{6S14f})}(\cdots)$ 
is defined by 
\begin{eqnarray}
\epsilon^{(\ref{6S14f})}_1(r,k,\vec{b},\epsilon;l_1,\cdots,l_{k-1}):=
\min\left\{
\epsilon_1^{(\ref{6O12d})}(\epsilon),
\epsilon_1^{(\ref{z0813})}
(r,k,\vec{L}(l_1,\cdots,l_{k-1}))
\right\}
\label{6O13}
\end{eqnarray}
where 
\begin{eqnarray}
\vec{L}(l_1,\cdots,l_{k-1}):=(L^{(\ref{6O12d})}(\epsilon,b'_i))_{i\in [k]}
\mbox{ and }
b'_i:=B_i(\vec{b},l_i+\cdots+l_{k-1}).
\label{6O14a}
\end{eqnarray}
\par
In this paragraph, we will define the function $\widetilde{m}(\cdot).$ 
Consider a sequence of integers $\vec{\ell}=(\ell_i)_{i\in [k-1]}$.
Theorem \ref{mt060721}
($r:=r,k:=k,\vec{b}:=\vec{b'}=(b'_i)_{i\in [k]},h:=h^\circ(b'_1,\cdots, b'_k),
\epsilon:=\epsilon_1^2=\epsilon_1^\circ(\vec{\ell})^2$) 
gives an $\widetilde{M}$ such that 
for any $\vec{\varphi}\in\Phi(\vec{\ell})$, 
there exist $m_1,\cdots,m_{k-1}\le \widetilde{M}$ for which 
\begin{eqnarray}
\Ex_{\vec{\psi}\in \Phi(m_1,\cdots,m_{k-1})}
[{\bf reg}_{k,h^\circ(b'_1,\cdots, b'_{k-1})}({\bf H}/\vec{\psi})]\le 
\epsilon_1^2,\label{6O12b}
\end{eqnarray}
 where ${\bf H}={\bf G}/\vec{\varphi}$ 
and $b'_i:=B_i(\vec{b},\ell_i+\cdots+\ell_{k-1})\ge c_i({\bf H})$ by (\ref{6O16d}).
Define 
$$\widetilde{m}(\vec{\ell}):= \widetilde{M}
.$$ 
\par
Next, we will define an integer $\tilde{\ell}$ as follows. 
Theorem \ref{6S19} 
($r:=r,k:=k,\vec{b}:=\vec{b},h:=1,
\epsilon:=0.1\epsilon, L(l_1,\cdots,l_{k-1})
:=(k-1)\widetilde{m}(l_1,\cdots,l_{k-1})$) 
 gives an integer $\tilde{\ell}$ such that 
(for any ${\bf G}$) there exist $\ell_1,\cdots,\ell_{k-1}\in [\tilde{\ell}]$ 
for which
\begin{eqnarray}
\Ex_{\vec{\varphi}\in\Phi(\ell_1,\cdots,\ell_{k-1})}
[{\bf reg}_{k,1,L(\ell_1,\cdots,\ell_{k-1})}
({\bf H}={\bf G}/\vec{\varphi})]\le 0.1\epsilon.
\label{6O12a}
\end{eqnarray}
\par
It suffices to show that these $\widetilde{\ell}$ and $\widetilde{m}(\cdots)$ 
satisfy the desired qualifications.
\par
For graph ${\bf G},$ there exist $\ell_1,\cdots,\ell_{k-1}\in [\widetilde{\ell}]$ 
satisfying (\ref{6O12a}).
Then we randomly pick a $\vec{\varphi}\in \Phi(\vec{\ell})$ with 
$\vec{\ell}=(\ell_1,\cdots,\ell_{k-1})$.
For this $\vec{\varphi}$, there exist $m_1,\cdots,m_{k-1}\in 
[\widetilde{m}(\vec{\ell})]$ satisfying (\ref{6O12b}).
Further we randomly pick a $\vec{\psi}\in \Phi(\vec{m})$ with 
$\vec{m}=(m_1,\cdots,m_{k-1})$.
By (\ref{6O12a}), for a random $\vec{\varphi}$, it holds 
with probability at least 0.9 that 
\begin{eqnarray}
\mbox{${\bf H}$ is $(\epsilon,k,1,(k-1)\widetilde{m}(\vec{\ell}))$-regular}.
\label{6O12c}
\end{eqnarray}
When (\ref{6O12c}) happens, 
since $m_1+\cdots+m_{k-1}\le (k-1)\widetilde{m}(\vec{\ell}),$
Lemma \ref{6O12d} with (\ref{6O12b}) 
gives 
positive-valued 
functions $\epsilon^{(\ref{6O12d})}(\cdot)$
 and $L^{(\ref{6O12d})}(\cdot,\cdot)$ such that if 
\begin{eqnarray}
0<\epsilon\le \epsilon^{(\ref{6O12d})}(k) 
\mbox{ and }0<\epsilon_1\le \epsilon_1^{(\ref{6O12d})}(\epsilon),
 \label{6O12e}
\end{eqnarray}
 then 
the $\vartheta$ probabilistically defined 
in Definition \ref{6S27} for 
 $\vec{L}(\ell_1,\cdots,\ell_{k-1})$ of (\ref{6O14a})
satisfies the inequality that 
$
\Ex_{\vec{\psi}\in \Phi(\vec{m})}
\Ex_{\mathfrak{\vartheta}}\Prob_{{\bf e}\in {\bf \Omega}_I}[
\vec{\vartheta}({\bf H}(\partial {\bf e}))
\mbox{ is $(\epsilon_1,\epsilon^{2/3},\epsilon^{1/3})$-ordinary }
]\ge 1-2^{|I|}\epsilon^{1/3}/c$ for all $I\in {\mathfrak{r}\choose [k]},$
which implies that 
\begin{eqnarray}
\Ex_{\vec{\psi}\in \Phi(\vec{m}
)}
\Ex_{\mathfrak{\vartheta}}\sum_{I\in {\mathfrak{r}\choose [k]}}
\Prob_{{\bf e}\in {\bf \Omega}_I}[
\vec{\vartheta}({\bf H}(\partial {\bf e}))
\mbox{ is not $(\epsilon_1,\epsilon^{2/3},\epsilon^{1/3})$-ordinary }
]\le 2^{r+k}\epsilon^{1/3}/c.
\label{6O12f}
\end{eqnarray}
Note that (\ref{6O12e}) is satisifed because of the assumption of the lemma 
and because of (\ref{6O13}).
Thus when (\ref{6O12c}) holds, for a random 
$\vec{\psi}\in\Phi(\vec{m})$, with probability at least $1-0.1-\epsilon_1\ge 0.89
$, 
it follows from (\ref{6O12f}) that 
\begin{eqnarray}
\Ex_{\mathfrak{\vartheta}}\sum_{I\in {\mathfrak{r}\choose [k]}}
\Prob_{{\bf e}\in {\bf \Omega}_I}[
\vec{\vartheta}({\bf H}(\partial {\bf e}))
\mbox{ is not $(\epsilon_1,\epsilon^{2/3},\epsilon^{1/3})$-ordinary }
]\le 10\cdot 2^{r+k}\epsilon^{1/3}/c
\label{6O12g}
\end{eqnarray}
and from $(\ref{6O12b})$ that
\begin{eqnarray}
\mbox{${\bf H}/\vec{\psi}$ is 
$(\epsilon_1,k,h^\circ(c_1({\bf H}),\cdots,c_{k-1}({\bf H})))$-regular}.
\label{6O12h}
\end{eqnarray}
By Lemma \ref{z0813} with (\ref{6O13}) and (\ref{6O12h}), we have that 
\begin{eqnarray}
\sum_{I\in {\mathfrak{r}\choose [k]}}
\sum_{\vec{\mathfrak{c}}\in {\rm TC}_I({\bf H})}
\Prob_{\vartheta}\left[
\vec{\vartheta}(\vec{\mathfrak{c}})
\in {\rm \overline{O}}^{\epsilon_1}{\rm TC}_I({\bf H}/
\vec{\psi})
\right]
&\le &
\sum_{I\in {\mathfrak{r}\choose [k]}}
\sum_{\vec{a}\in {\mathcal A}_I}
\Prob_{\mathfrak{d}}\left[
\vec{\mathfrak{d}}(\vec{a})
\in {\rm \overline{O}}^{\epsilon_1}{\rm TC}_I({\bf H}/
\vec{\psi})
\right]
\nonumber\\
&\le &
 0.01.\label{6O14}
\end{eqnarray}
Thus by (\ref{6O12g}) and (\ref{6O14}), for a random process 
of $\vartheta$, with probability at least $1-0.01-0.001\ge 0.9,$
the desired properties (i) and (ii) hold simultaneously.
\par
It easily follows from the definition of $\vec{\vartheta}$ that if 
${\mathfrak{c}}=(\mathfrak{c}_J)_{J\subsetneq I}
\in \partial {\rm C}_I({\bf H})$ and 
$\vec{\vartheta}(\mathfrak{c})\not=\mathfrak{0}$  then there exists a color 
$\mathfrak{c}_I\in {\rm C}_I({\bf H})$ such that 
$\vec{\vartheta}((\mathfrak{c}_J)_{J\subset I})\not=\mathfrak{0}$.
Thus property (ii) implies  (iii).
It completes the proof of Lemma \ref{6S14f}.
\lemmaqed
\section{Body Part of the Proof of Main Theorem}
We will prove Theorem \ref{6S05} by using Lemmas \ref{6O12d} and \ref{6S14f}.
\bigskip\\
{\bf Proof of Theorem \ref{6S05}: }
Let $r,k,\vec{b},\varepsilon$ be given as in the theorem.
(Without loss of generality, $\vec{b}=(1,\cdots,1,b_k)$, though we will not use this.)
Let $0<\epsilon\le \epsilon^{(\ref{6O12d})}(k),$ and $h^\circ(\cdot)
$ be a function, which will be defined later at (\ref{6O16}) and 
 (\ref{6S14h}), respectively.
Let $0<\epsilon_1^\circ(\cdot)=\epsilon_1^{(\ref{6S14f})}(\cdot)$ be the 
function which decreases fast enough in Lemma \ref{6S14f}.
By Lemma \ref{6S14f} with $r,k,\vec{b},\epsilon,\epsilon_1^\circ, 
h^\circ$ and with ${\bf G}$,
there exist an integer $\widetilde{\ell}$ and a 
function $\widetilde{m},$ which are independent from ${\bf G},$ 
together with 
$\vec{\varphi}\in \Phi(\vec{\ell})$ and $\vec{\psi}\in \Phi(\vec{m})$  
 for some $\vec{\ell}\in [\widetilde{\ell}]^{k-1}$ 
and $\vec{m}\in [\widetilde{m}(\vec{\ell})]^{k-1}$ 
such that ${\bf H}/\vec{\psi}$ is 
$\brkt{\epsilon_1:=\epsilon_1^\circ(\vec{\ell}),k,h^\circ(
c_1({\bf H}),\cdots,c_{k-1}({\bf H}))
}$-regular, where
 ${\bf H}:={\bf G}/\vec{\varphi}.$ 
Furthermore there exist a map $\vec{\vartheta}$ which satisfies 
properties (i)-(iii) of Lemma \ref{6S14f} simultaneously.
\bigskip\\
{\bf [Modification of ${\bf G}$]}\quad 
By conducting the steps $S_1,\cdots,S_k$, which will be defined below,
we will redefine the face-colors 
${\bf H}({\bf e})$ for edges ${\bf e}\in {\bf \Omega}_I, 
I\in {\mathfrak{r}\choose [k]}$. We will denote the new colored hypergraph by ${\bf H'}$, 
instead of ${\bf H}$.(We will see ${\rm C}_I({\bf H})={\rm C}_I({\bf H'})$ since 
we will not add any new color, and will not remove any unused color 
from the color class, either. 
We always use simbol ${\bf H}$ 
for the old one.)
\par
(Step $S_s$) 
Assume that ${\bf H'}({\bf e'})$ has been defined for all ${\bf e'}\in 
\bigcup_{J\in {\mathfrak{r}\choose [s-1]}}{\bf \Omega}_J$ so that 
 $\vec{\vartheta}({\bf H'}\ang{\bf e'})\in 
{\rm O}^{\epsilon_1}{\rm TC}_J({\bf H}/\vec{\psi})$.
Let 
$I\in {\mathfrak{r}\choose s}$ and ${\bf e}\in {\bf \Omega}_I$. 
Write $\mathfrak{c}_I:={\bf H}({\bf e})$ and 
${\mathfrak{c}}=({\mathfrak c}_J)_{J\subsetneq I}:={\bf H'}(\partial {\bf e})
\in \partial {\rm C}_I({\bf H})$.
By the assumption for $s-1$, $\vec{\vartheta}(\mathfrak{c})\in 
{\rm O}^{\epsilon_1}
\partial {\rm C}_I({\bf H}/\vec{\psi})$.
Our purpose of this step is to define face-color ${\bf H'}({\bf e})\in 
{\rm C}_I({\bf H}).$ 
\\
(Case $S_s-1$) Suppose that $\vec{\vartheta}(\mathfrak{c})$ is $(\epsilon_1,
\epsilon^{2/3},\epsilon^{1/3})$-ordinary and that 
${\bf d}_{{\bf H}}(\mathfrak{c}_I|\mathfrak{c})\ge {\epsilon^{1/3}\over 
|{\rm C}_I({\bf H})|}$.
Define 
$${\bf H'}({\bf e}):=\mathfrak{c}_I.$$ 
(Case $S_s-1'$) Suppose that $\vec{\vartheta}(\mathfrak{c})$ is $(\epsilon_1,
\epsilon^{2/3},\epsilon^{1/3})$-ordinary and that 
${\bf d}_{{\bf H}}(\mathfrak{c}_I|\mathfrak{c})< {\epsilon^{1/3}\over 
|{\rm C}_I({\bf H})|}$.
Note that there exists such a color  $\mathfrak{c}_I^\circ \in {\rm C}_I({\bf H})$ 
such that ${\bf d}_{{\bf H}}(\mathfrak{c}_I^\circ
|\mathfrak{c})\ge {1\over 
|{\rm C}_I({\bf H})|}>{\epsilon^{1/3}\over 
|{\rm C}_I({\bf H})|}$. Fix such a color and define 
$${\bf H'}({\bf e}):=\mathfrak{c}_I^\circ.$$ 
(Case $S_s-2$) Suppose that $\vec{\vartheta}(\mathfrak{c})$ is not $(\epsilon_1,
\epsilon^{2/3},\epsilon^{1/3})$-ordinary.
(This case does not occur when $s=1.$) 
Then by (iii) of Lemma \ref{6S14f}, 
there exists a color $\mathfrak{c}_I^\circ\in {\rm C}_I({\bf H})$ such that 
$\vec{\vartheta}((\mathfrak{c}_J^\circ)_{J\subset I})\in 
{\rm O}^{\epsilon_1}{\rm TC}_I({\bf H}/\vec{\psi})$ 
where $\mathfrak{c}^\circ_J:=\mathfrak{c}_J$ for all $J\subsetneq I$.
Fixing one, we define 
$${\bf H'}({\bf e}):=\mathfrak{c}_I^\circ.$$ 
\par
By Lemma \ref{6S14f} (ii) with the definition of $(*,*,\epsilon^{1/3})$-ordinarity 
(Lemma \ref{6O12d} (iii)), 
we see that $
\vec{\vartheta}({\bf H'}\ang{\bf e})\in 
{\rm O}^{\epsilon_1}{\rm TC}_I({\bf H}/\vec{\psi})$ 
for any of the three cases.
\bigskip\\
{\bf [Estimating the edit size]}
We define value ${\rm Ordinariness}({\bf e})\in [0,|I|]$ 
for  ${\bf e}\in {\bf \Omega}_I$ 
by the largest integer $s\ge 0$ such that 
for any $J\subset I$ with $|J|\le s,$ 
${\bf H'}({\bf e}|_J)$ was defined by (Case $S_{|J|}-1$).
Note that if ${\rm Ordinariness}({\bf e})=|I|$ then 
${\bf H'}\ang{\bf e}={\bf H}\ang{\bf e}$ by Lemma \ref{6O12d} (iii).
For any $I\in {\mathfrak{r}\choose k},$
we have that 
\begin{eqnarray}
\Prob_{{\bf e}\in {\bf \Omega}_I}[{\bf H'}({\bf e})\not={\bf G}({\bf e})
]
&\le &
\Prob_{{\bf e}\in {\bf \Omega}_I}[{\bf H'}\ang{\bf e}\not={\bf H}\ang{\bf e}
]
\nonumber\\
&\le &
\Prob_{{\bf e}\in {\bf \Omega}_I}[
{\rm Ordinariness}({\bf e})<k
]
\nonumber\\
&\le &
\sum_{J\subset I}
\Prob_{{\bf e}\in {\bf \Omega}_I}[
{\rm Ordinariness}({\bf e}|_J)=|J|-1
]
\nonumber\\
&\le &\sum_{J\subset I}
\Prob_{{\bf e}\in {\bf \Omega}_J}[
{\rm Ordinariness}({\bf e})\not=|J|
|
{\rm Ordinariness}({\bf e})\ge |J|-1
]
\nonumber\\
&\le &
\sum_{J\subset I}
\Prob_{{\bf e}\in {\bf \Omega}_J}[
{\bf H}(\partial {\bf e})\mbox{ is not
 $(\epsilon_1,\epsilon^{2/3},\epsilon^{1/3})$-ordinary}
|
{\rm Ordinariness}({\bf e})\ge |J|-1
]
\nonumber\\
&&
+
\sum_{J\subset I}
\Prob_{{\bf e}\in {\bf \Omega}_J}[
{\bf d}_{{\bf H}}({\bf H}\ang{\bf e})<{\epsilon^{1/3}\over 
|{\rm C}_J({\bf H})|}
|
{\rm Ordinariness}({\bf e})\ge |J|-1
]
\nonumber\\
&\le &\sum_{J\subset I}O_{r,k}(\epsilon^{1/3}) +
\sum_{J\subset I}\epsilon^{1/3} \quad\quad 
(\because \mbox{ Lemma } \ref{6S14f}\ (i))
\nonumber\\
&= & O_{r,k}(\epsilon^{1/3})
\nonumber\\
&\le &\varepsilon,\label{6O16}
\end{eqnarray}
when $\epsilon$ is small enough for $r,k,\varepsilon.$
\\
{\bf [Choosing a target forbidden graph]} 
When $b'_i, i\in [k-1],$ are integers 
and when $g_I:\prod_{J\subsetneq I}[b'_{|J|}]
\to 2^{[b_k]}\setminus \{\emptyset\}
, I\in {\mathfrak{r}\choose k},$ are maps,
we say that an $F'\in {\mathcal F}_h$ is 
{\bf $((b'_i)_i, (g_I)_I)$-colorable} if and only if there exists a $k$-bound 
$(b'_1,\cdots,b'_{k-1},b_k)$-colored 
simplicial-complex $S\in {\mathcal S}_{\mathfrak{r},k,h}$ 
on the same vertex sets (as $F'$) such that 
 $F'(e)=S(e)\in g_I(S(\partial e))$ for all $e\in \mathbb{V}_I(F'), 
I\in {\mathfrak{r}\choose k}.$ 
Given integers $M_1,\cdots,M_{k-1}$, we define 
\begin{eqnarray}
h^\circ(M_1,\cdots,M_{k-1}):=h_0\label{6S14h}
\end{eqnarray}
 to be the smallest value $h_0\ge 0$ such that 
for any $(b'_1,\cdots,b'_{k-1})$ with $b'_i\le M_i,\forall i\in [k-1],$ 
and for any $(g_I)_{I\in {\mathfrak{r}\choose k}}$, 
at least one of the following two holds:
\begin{itemize}
\item[(a)] There does not exist a $((b'_i)_i, (g_I)_I)$-colorable graph 
$F'\in {\mathcal F}$.
\item[(b)]
There exists a $((b'_i)_i, (g_I)_I)$-colorable graph 
$F'\in {\mathcal F}_{h}$ with $h\le h_0$ (or $h=h_0$ without loss of generality, 
by adding extra invisible edges).
\end{itemize}
\par
Assume that (i) of the theorem does not hold. Then 
there exist an $h\ge 1$ and an $F\in {\mathcal F}_h$ such that 
\begin{eqnarray*}
\Prob_{\phi\in\Phi(h)}[{\bf H'}(\phi(e))=F(e)\forall e\in {\mathbb V}(F)]>0.
\end{eqnarray*}
By the image of a map $\phi\in\Phi(h)$ with the above property, we can construct an 
$S\in {\mathcal S}_{r,k,h,{\bf H'}}$ which shows the $((b'_i)_i, (g_I)_I)$-colorability of 
$F$ where 
\begin{eqnarray*}
b'_i&:=&c_i({\bf H'})= c_i({\bf H})
,
\\
g_I=g_I(\mathfrak{c})&:=&
\{\mathfrak{c}_I\in {\rm C}_I({\bf H'})= {\rm C}_I({\bf G})
|
 {\bf d}_{\bf H'}(\mathfrak{c}_I| \mathfrak{c})>0
\}
\quad \forall \mathfrak{c}\in \partial {\rm C}_I({\bf H'})
\end{eqnarray*}
 (under some map
from ${\rm C}_{I'}({\bf H'})= {\rm C}_{I'}({\bf H})$ 
to $[b'_i]$).
By the definition (\ref{6S14h}) of $h^\circ$ and by the existence of colorable $F
\in {\mathcal F},$ for the above pair $((b'_i)_i,(g_I)_I))$,
the item (a) does not happen, and then 
there exists a $((b'_i)_i,(g_I)_I)$-colorable $F^*\in {\mathcal F}_{h_0}$ 
where 
\begin{eqnarray*}
h_0:=h^\circ(b'_1,\cdots,b'_{k-1}
),
\end{eqnarray*}
which is smaller than a constant depending only on 
$r,k,\vec{b},\varepsilon$ and ${\mathcal F}$ since 
$h^\circ (\cdot)$ is monotone without loss of generality and 
\begin{eqnarray}
b_i'
\stackrel{(\ref{6O16d})}{\le} 
 B_i(\vec{b},(k-i)\widetilde{\ell})
.
\label{6O16e}
\end{eqnarray}
Let $S^*\in {\mathcal S}_{r,k,h_0,{\bf H'}}$ be the simplicial-complex guaranteeing the 
colorability of $F^*$.
\bigskip\\
{\bf [Finding many copies]} 
We will show that there exist many copies of $F^*$ in ${\bf G}$.
For this purpose, we define $S^{**}\in {\mathcal S}_{r,k,h_0,{\bf H}/\vec{\psi}}$ 
from $S^*\in {\mathcal S}_{r,k,h_0,{\bf H'}}$ 
by replacing $S^*(e)$ by $S^{**}(e):={\vartheta}(S^*\ang{e})$ for each $e\in 
\mathbb{V}(S^*)=\mathbb{V}(S^{**})$.
\par
By the definition of $\vartheta(\cdot)$, 
if $|I|=k
$ and 
$\vec{\mathfrak{c}}=(\mathfrak{c}_J)_{J\subset I}
\in {\rm TC}_I({\bf H'})$ 
then $\vartheta(\vec{\mathfrak{c}})\in \{\mathfrak{o}, \mathfrak{c}_I\}$ since 
${\rm C}_I({\bf H'})={\rm C}_I({\bf H}/\vec{\psi})$.
Therefore by our definiton of ${\bf H'}$, 
if ${\bf e}\in {\bf \Omega}_I$ with $I\in {\mathfrak{r}\choose k}$ 
then $\mathfrak{o}\not=\vartheta({\bf H'}\ang{\bf e})={\bf H'}({\bf e}).$
Using this fact, it is easily seen that not only $S^*$ but also 
$S^{**}$ is a simplicial-complex guaranteeing the $((b'_i)_i,(g_I)_I)$-colorability 
of $F^*$ by identifying $S^*(e)$ as $S^{**}(e)=\vartheta(S^*\ang{e})$ 
(in the domain of $g_I$) 
for each $e\in 
\mathbb{V}_I(S^*), I\in {\mathfrak{r}\choose [k]}.$
(To see this, for all 
$e\in \mathbb{V}_k(F^*)$, 
observe that $F^*(e)=S^*(e)=\vartheta(S^*\ang{e})
=S^{**}(e)$ and that 
$S^*(e)\in g_I(S^{*}(\partial e))
\stackrel{\rm identify}{=}
g_I(\vec{\vartheta}(S^{*}(\partial e)))=g_I(S^{**}(\partial e)).
$)
\par
By Lemma \ref{6S14f} (ii) with the definition of $(*,*,\epsilon^{1/3})$-ordinarity 
(Lemma \ref{6O12d} (iii)), 
we have that $\vec{\vartheta}(S^*\ang{e})\in 
{\rm O}^{\epsilon_1}
{\rm TC}_I({\bf H}/\vec{\psi})$ for 
all $e\in \mathbb{V}_I(S^*), I\in {\mathfrak{r}\choose [k]}.$ 
Hence it follows 
from 
\begin{eqnarray}
c_i({\bf H}/\vec{\psi})
\stackrel{(\ref{6O16d})}{\le} 
 B_i((b'_j)_j,(k-i)\widetilde{m}(\widetilde{\ell},\cdots,\widetilde{\ell}))
\label{.101}
\end{eqnarray}
that 
\begin{eqnarray*}
 \Prob_{\phi\in \Phi(h_0)}[{\bf G}(\phi(e))=
F^*(e)
\forall e\in \mathbb{V}(F^*)]
&\ge&
 \Prob_{\phi\in \Phi(h_0)}[{\bf H}/\vec{\psi}(\phi(e))=
S^{**}(e)\forall e\in \mathbb{V}(S^{**})]\\
&=&\prod_{e\in \mathbb{V}(S^{**})}({\bf d}_{{\bf H}/\vec{\psi}}(S^{**}\ang{e})
\dot{\pm}\delta(S^{**}\ang{e}))
\\
&=&\prod_{I\in {\mathfrak{r}\choose [k]}}
\prod_{e\in \mathbb{V}_I(S^*)}({\bf d}_{{\bf H}/\vec{\psi}}(
\vec{\vartheta}(S^*\ang{e}))
\dot{\pm}\delta(\vec{\vartheta}(S^*\ang{e}))
\\
&\ge &\prod_{I\in {\mathfrak{r}\choose [k]}}
\prod_{e\in \mathbb{V}_I(S^*)}
{
\epsilon_1^{1/3}-\epsilon_1^{2/3}
\over |{\rm C}_I({\bf H}/\vec{\psi})|
}
\\
&\stackrel{(\ref{.101})}{\ge} &\prod_{i=1}^k
\brkt{0.9
(\epsilon_1^{(\ref{6S14f})}(\widetilde{\ell},\cdots,\widetilde{\ell})
)^{1/3}
\over B_i((b'_j)_j,(k-i)\widetilde{m}(\widetilde{\ell},\cdots,\widetilde{\ell}))
}^{{r\choose i}h_0^i},
\end{eqnarray*}
which is larger than 
a positive real depending only on $r,k,\vec{b},\varepsilon$ and ${\mathcal F}$ 
by (\ref{6O16e}). 
In the last inequality, we used the fact that function 
$\epsilon_1^{(\ref{6S14f})}\le 0.1^{3/2}$ is monotone without loss of generality.
It completes the proof of Theorem \ref{6S05}.
\qed


\begin{thebibliography}{99}
\bibitem{AFKS}N.~Alon, E.~Fischer, M.~Krivelevich and M.~Szegedy,
Efficient testing of large graphs,
{\em Proc. 40th FOCS, New York, NY, IEEE} (1999), 656-666. Also: {\em Combinatorica}
{\bf 20} (2000), 451-476.
\bibitem{AS}
N. Alon and A. Shapira,
Homomorphisms in graph property testing - a survey,
{\em Topics in Discrete Mathematics,} 28 pages, to appear. 
\bibitem{AS05}
N. Alon and A. Shapira,
Every monotone graph property is testable,
{\em Proc. of STOC 2005}, 128-137.
(Also, {\em SICOMP (Special Issue of STOC'05)}, to appear.)
\bibitem{ARS}C.~Avart, V.~R\"{o}dl and M.~Schacht, 
Every monotone 3-graph property testable, 21 pages, preprint (2006.2).
\bibitem{BLR}M.~Blum, M.~Luby and R.~Rubinfeld, 
Self-testing/correcting with applications to numerical problems, 
{\em Journal of Computer and System Sciences} {\bf 47} (1993), 549-595
(a preliminary version appeared in Proc. 22nd STOC, 1990)
\bibitem{CRST}V.~Chv\'atal, V.~R\"{o}dl, E.~Szemer\'edi, and
W.T.Trotter, Jr., The Ramsey number of a graph with bounded maximum
degree, {\it J. Combin. Theory, Ser. B}, {\bf 34} (1983), 239-243.
\bibitem{CFKO2}O.~Cooley, N.~Fountoulakis, D.~K\"{o}hn, and D.~Osthus,
Embeddings and Ramsey numbers of sparse k-uniform hypergraphs,  
preprint, 
{\tt 
arXiv:math/0612351v1 [math.CO]}.
\bibitem{Fis}E.~Fischer, The art of uninformed decisions: A primer to property testing, 
{\em The Bulletin of the European Association for Theoretical Computer Science}
 {\bf 75} (2001), 97-126.
\bibitem{FR02}P. Frankl and V. R\"{o}dl, Extremal problems on set systems, 
{\it Random Structures and Algorithms}, {\bf 20}(2), 131-164 (2002).
\bibitem{FK78}H. Furstenberg and Y. Katszenelson, 
An ergodic Szemer\'edi theorem for 
commuting transformations, {\it J. Analyse Math.} {\bf 34} (1978), 275-291.
\bibitem{GGR}O.~Goldreich, S.~Goldwasser and D.~Ron, 
Property testing and its connection to learning and approximation, 
{\em Journal of the ACM} {\bf 45} (1998), 653-750 (a preliminary version appeared in 
Proc. 37th FOCS, 1996)
\bibitem{GT}O.~Goldreich and L.~Trevisan, Three theorems
regarding testing graph properties, {\em Random Structures and
Algorithms}, {\bf 23} (2003), 23-57.
\bibitem{G}W.T.~Gowers, Hypergraph regularity and the multidimensional Szemer\'edi 
	theorem, 42 pages, preprint (2005.4, 2nd ver.)
\bibitem{I06}Y.~Ishigami, A simple regularization of hypergraphs, 13 pages, 
preprint, 
{\tt arXiv:math/0612838v1} [math.CO].
\bibitem{I06lr}Y. Ishigami, Linear Ramsey numbers for bounded-degree hypergraphs, 
preprint, 
{\tt arXiv:math/0612601v1} [math.CO].
\bibitem{I06m}Y.~Ishigami, Removal lemma for infinitely-many forbidden hypergraphs 
and property testing, preprint, {\tt arXiv:math/0612669v1} [math.CO].
\bibitem{KSSS02} J. Koml\'os, A. Shokoufandeh, M. Simonovits, and 
E.Szemer\'edi, The regularity lemma and its applications in graph theory, 
{\it Theoretical Aspects of Computer Science.} (Edited by G.B.Khosrovshahi et al.) 
 Lecture Notes in Computer Science Vol. 2292 (2002), 84-112.
\bibitem{KR}A.V.Kostochka and V.R\"{o}dl, On Ramsey numbers of 
uniform hypergraphs with given maximum degree, 
{\em Journal of Combinatorial Theory}, {\bf A 113} (2006) 1555-1564.
\bibitem{LS04}L.~Lov\'asz andd B.~Szegedy, Limits of dense graph sequences, 
{\em Journal of Combinatorial Theory}, {\bf B 96} (2006) 933-957.
\bibitem{LS05}L.~Lov\'asz andd B.~Szegedy, Graph limits and testing hereditary 
graph properties, {\em Tech. Report MSR-TR-2005-110}, Microsoft Research, 2005.
\bibitem{NRS}B.~Nagle, V.~R\"{o}dl and M.~Schacht, 
The counting lemma for regular $k$-uniform hypergraphs, 
{\em Random Structures and Algorithms}, {\bf 28} (2006), no.2, 113-179.
\bibitem{RSgene}V.R\"{o}dl and M.Schacht, Generalizations of the removal lemma, 25 pages, 
preprint (2006).
\bibitem{RSgene2}V.R\"{o}dl and M.Schacht, Property testing in hypergraphs and the removal lemma, 
(extended abstract), preprint (2006.11).
\bibitem{RS}V.R\"{o}dl and M.Schacht, Regular partitions of hypergraphs, 
{\em Combinatorics, Probability \& Computing}, to appear (preprint 50 pages, 2006.5).
\bibitem{RSk04}V.~R\"{o}dl and J.~Skokan, Regularity lemma for 
$k$-uniform hypergraphs, {\em Random Structures and Algorithms} {\bf 25} (2004) (1),
1-42.
\bibitem{Ron}D.~Ron, Property testing (a tutorial), In: {\em Handbook of 
Randomized Computing} (S.~Rajasekaran, P.M.~Pardalos, J.H.~Reif 
and J.D.P.~Rolin eds), Kluwer Press (2001)
\bibitem{RuSu}R.~Rubinfeld and M.~Sudan, Robust characterizations of 
polynomials with applications to program testing, {\em SIAM J.~Comput.} 
{\bf 25} (1996), no.2, 252-271.
\bibitem{RuSz}I.Z.~Ruzsa and E.~Szemer\'edi, Triple systems with 
no six points carrying three triangles, 
{\em Combinatorics (Proc. Fifth Hungarian Colloq., Keszthely, 1976)}, Vol.II,
Colloq.~Math.~Soc.~J\'anos~Bolyai, vol.~18, North-Holland, Amsterdam, 1978, pp.939-945.
\bibitem{So03}J.~Solymosi, Note on a generalization of Roth's theorem, 
{\em Discrete and Computational Geometry}, 825-827, {\em Algorithms Combin.} 
{\bf 25}, Springer, Berlin 2003.
\bibitem{So04}J.~Solymosi, A note on a question of Erd\H{o}s and Graham, 
{\em Combin. Probab. Comput.} {\bf 13} (2004), 263-267.
\bibitem{Sz69}E.~Szemer\'edi, On sets of integers containing no four elements 
in arithmetic progression, 
{\em Acta Math. Acad. Sci. Hungar.} {\bf 20} (1969), 89-104.
\bibitem{Sz75}E.~Szemer\'edi, On sets of integers containing no $k$ elements 
in arithmetic progression, 
{\it Acta Arithmetica} {\bf 27} (1975), 199-245. [Collection of articles in memory of
 Juri\u{i} Vladimirovi\u{c} Linnik.]
\bibitem{Sz}
	E. Szemer\'edi,
	{Regular partitions of graphs}
	in {\it Probl\`emes combinatoires et th\'eorie des graphes},
	Orsay 1976,
	J.-C. Bermond, J.-C. Fournier, M. Las Vergnas, D. Sotteau, eds.,
	Colloq. Internat. CNRS 260,
	Paris, 1978, 399--401.
\bibitem{Tao06}T.~Tao, A variant of the hypergraph removal lemma, 
{\em J.~Combin.~Theory} {\bf A 113} (2006), no.7, 1257-1280.
\bibitem{Tao06a}T.~Tao, 
The dichotomy between structure and randomness, arithmetic progressions, and the primes, 
(ICM2006 lecture) (preprint 27 pages, 2005.12, ver.2)
\bibitem{Tao06b}T.~Tao, 
A correspondence principle between (hyper)graph theory and probability 
theory, and the (hyper)graph removal lemma, 40 pages, preprint (2006.2)
\bibitem{TV}T.~Tao and V.H.~Vu, {\em Additive Combinatorics}, Cambridge University Press, 
(2006) 512 pages.
\end{thebibliography}
\end{document}